\documentclass[a4paper,12pt,oneside,reqno]{amsart}
\usepackage[headinclude,DIV13]{typearea}
\areaset{15.1cm}{25.0cm}
\usepackage{txfonts,amssymb,amsmath,amsthm,bbm,physics}
\usepackage{cases}
\usepackage{subfiles}

\usepackage[latin1] {inputenc}
\usepackage{graphicx, psfrag}
\usepackage{setspace}
\usepackage{xcolor}
\usepackage{subcaption}
\usepackage{enumerate}

\graphicspath{{Figures/}}

\usepackage{appendix}

\usepackage{hyperref}

\newtheorem{theorem}{Theorem}[section]
\newtheorem{lemma}[theorem]{Lemma}

\newtheorem{corollary}[theorem]{Corollary}

\DeclareMathOperator*{\argmax}{arg\,max}

\newcommand{\R}{\mathbb{R}}
\newcommand{\Prob}[1]{\mathbb{P}\left(#1\right)}

\renewcommand{\l}[2]{\ensuremath{I_{#1}\left(#2\right)}}
\renewcommand{\ll}[3]{\ensuremath{I_{#1}\left(#2;\,#3\right)}}

\definecolor{bbm}{RGB}{51,153,0}
\definecolor{above}{RGB}{128,0,128}
\definecolor{below}{RGB}{102,0,204}
\definecolor{cascade}{RGB}{204,0,0}
\definecolor{iid}{RGB}{153,51,0}

\theoremstyle{remark}

\def\paragraph#1{\noindent \textbf{#1}}

\numberwithin{equation}{section}

\def\dd{\mathtt{d}}

\def\<{\langle}
\def\>{\rangle}

\def\a{\alpha}

\def \g {{\gamma}}

\def \a {{\alpha}}

\def \1{\mathbbm{1}}

\def\eee{{\mathrm e}}

\def \log{\ln}

\begin{document}


 \title[LDP for BBM ]{Refined large deviation principle for branching Brownian motion conditioned to have a low maximum}
\author[Y. Bai]{Yanjia Bai}
 \address{Y. Bai\\Institut f\"ur Angewandte Mathematik\\
Rheinische Friedrich-Wilhelms-Universität\\ Endenicher Allee 60\\ 53115 Bonn, Germany }
\email{bai@iam.uni-bonn.de}
\author[L. Hartung]{Lisa Hartung}
 \address{L. Hartung\\
  Institut für Mathematik \\Johannes Gutenberg-Universität Mainz\\
Staudingerweg 9,
55099 Mainz, Germany}
\email{lhartung@uni-mainz.de}

\date{\today}

 \begin{abstract}  
Conditioning a branching Brownian motion to have an atypically low maximum leads to a suppression of the branching mechanism. In this note, we consider a branching Brownian motion conditioned to have a maximum below $\sqrt{2}\alpha t$ ($\a<1$). We study the precise effects of an early/late first branching time and a low/high first branching location under this condition. We do so by imposing additional constraints on the first branching time and location. We obtain large deviation estimates, as well as the optimal first branching time and location given the additional constraints.

\end{abstract}

\thanks{
The work of Y. B. is partially funded by the Deutsche Forschungsgemeinschaft (DFG, German Research Foundation) - Project-ID 211504053 -
SFB 1060 and German's Excellence Strategy - GZ 2047/1, Project-ID 390685813 -"Hausdorff Center for Mathematics" at
Bonn University. The research of L. H. is supported in part by the Deutsche Forschungsgemeinschaft
(DFG, German Research Foundation) through Project-ID 233630050 -TRR 146,
Project-ID 443891315 within SPP 2265 and and Project-ID 446173099. }

\subjclass[2010]{60J80, 60G70, 82B44} \keywords{branching Brownian motion, large deviation principle, extreme values} 

 \maketitle


	\section{Introduction  }
	
%
	
	In this note, we contribute to a deeper understanding of how spatial branching processes or log-correlated Gaussian processes realize unlikely events, such as having a low maximum. For a continuous time branching processes, the particles can refrain from branching with only exponential costs. As a result, when a continuous time branching process is conditioned to behave atypically, there is an interesting interplay of three factors: the suppression of branching, the first branching time, and the first branching location.  
	
	A branching Brownian motion (BBM) $X$ can be constructed as follows. Starting from the origin at time $0$, one particle performs an one dimensional standard Brownian motion. After an exponentially distributed time with parameter one, the initial particle splits into two particles. From this branching location, the new particles follow independent Brownian paths and are subject to the same splitting rule. We denote the number of particles at time $t$ by $n(t)$ and the particle positions by $\{X_k(t):1\leq k\leq n(t) \}$. 

	We study the aforementioned three-factor interplay for BBM by asking the maximum to be a linear order below its typical value, while constraining the first branching time and location. When putting repulsive constraints, the branching mechanism is suppressed and the BBM has fewer particles. This seems to be a universal phenomenon. The repulsive constraints could come through spatial inhomogeneities \cite{E08,EH03,EK04, OE19}, or a direct repulsive interaction through the center of mass (as conjectured in \cite{E10}), or a polymer-type change of measure \cite{BH20}, or a requirement of an unusual maximum \cite{chen2020branching, derrida2016large}.

	The extreme values of branching Brownian motion have been studied extensively over the last 40 year (see, e.g. \cite{ABBS, ABK_G,ABK_E,B_M,B_C,CHL17}), and it is well known that the order of the maximum is given by
	\begin{equation}\label{int.1}
	m_t:=\sqrt{2}t-\frac{3}{2\sqrt{2}}\log t.
	\end{equation}	
	Let $X_{max}(t):=\max_{k\leq n(t)} X_k(t)$ denote the position of the maximal particle at time $t$. The probability of BBM having a maximum smaller than $\sqrt{2}\a t$, $\a<1$, was first analyzed by Derrida and Shi \cite{derrida2016large}. They showed that 
		\begin{equation}\label{ds}
			\lim_{t\rightarrow\infty}\frac{1}{t}\log(\Prob{X_{max}(t)\leq \sqrt{2}\alpha t})= -\psi(\alpha),
			\quad
			\psi(\alpha):=
			\begin{cases}
				2\rho(1-\alpha),&\text{if }\alpha\in[-\rho,1),\\
				1+\alpha^2,&\text{if }\alpha\in(-\infty,-\rho],
			\end{cases}
		\end{equation}
	where $\rho:=\sqrt{2}-1$. Moreover, they pointed out that if $\tau$ denotes the first branching time and $y$ denotes the particle location at time $\tau$, then the optimal choices are
	\begin{equation}\label{ds_optimal_tau}
	\tau= \left(\frac{1-\alpha}{\sqrt{2}}\wedge 1\right)t+o(t)
	\end{equation}
	and
	\begin{equation}\label{ds_optimal_y}
	y=\sqrt{2}\alpha t-\sqrt{2}(t-\tau)+o(t)=\begin{cases}
	-\rho(1-\alpha)t+o(t),&\text{ if }\alpha\in[-\rho,1),\\
	\sqrt{2}\alpha t+o(t),&\text{ if }\alpha\in(-\infty,-\rho].
	\end{cases}
	\end{equation}
	These results were further refined by Chen, He, and Mallein \cite{chen2020branching}. They gave precise constant and polynomial prefactors of the probability in \eqref{ds}, and proved the convergence in distribution of the first branching time and first branching location conditioned on the BBM having a low maximum. Moreover, they showed convergence of the extremal process under this conditioning.
	
	The particular case $\alpha=0$, which restricts the field to be below zero, has also received some attention for models related to BBM, such as the $2d$ discrete Gaussian free field and the binary branching random walk.   In \cite{BDG01}, Bolthausen, Deuschel, and Giacomin analyzed the $2d$ discrete Gaussian free field conditioned to be below zero away from the boundary.  In \cite{Roy18}, Roy studied a binary branching random walk conditioned to be below zero and gave bounds on the height of a typical vertex under the conditioning.   Both models fall into the same universality class as BBM on the level of extremes.

	Although this paper focuses on the BBM whose maximum is unusually low, it is worth noting that the probability of BBM to have an unusually high maximum  $u(t,x(t)+\sqrt{2} t)\equiv \Prob{X_{max}(t)\geq x(t)+\sqrt{2} t}$ has been studied (see, e.g. \cite{ABK_G,B_C,chauvin88,chauvin90,derrida2016_physics,harris1999,LS}). $u(t,x)$ solves the Kolmogorov-Petrovsky-Piscounov or Fisher (F-KPP) equation \cite{fisher,kpp}  (see \cite{Ikeda1,Ikeda2,Ikeda3, McKean}). 
	The asymptotics of $u(t,x(t)+\sqrt{2}t)$ have been obtained for different ranges of $x(t)$ based on Bramson's analysis \cite{B_C}. 
	If $x(t)=o(t)$, then (see Bovier and Hartung \cite[Proposition 2.1]{BH1to6})
	\begin{equation}
		\lim_{t\rightarrow\infty}\frac{t^{\frac{3}{2}}}{\frac{3}{2\sqrt{2}}\log(t)}\eee^{\sqrt{2}x(t)+\frac{x(t)^2}{2t}}u\left(t,x(t)+\sqrt{2}t \right)=C,
	\end{equation}
	and if $x(t)=at+o(t), a>0$, then (see \cite[Proposition 3.1]{BovHar13})
	\begin{equation}\label{u_3}
		\lim_{t\rightarrow\infty}t^{\frac{1}{2}}\eee^{\sqrt{2}x(t)+\frac{x(t)^2}{2t}}u\left(t,x(t)+\sqrt{2}t \right)=C(a),
	\end{equation}
	where $C$ and $C(a)$ are strictly positive constants.
	
	\subsection{Main results}

	 Intuitively, three types of behaviors of  BBM may lead to it having a low maximum: the initial particle branches at a late time, the initial particle travels to a low position before branching, or the two  independent BBMs starting from the first branching position both have low maxima.  In this paper, we are interested in the interplay of these effects and aim at quantifying the exact scale of the decay in the large deviation estimates, with restrictions on the first branching time and location.
	 
Set $\tau:=\inf\{0\leq s\leq t,\, n(s)>1 \}$ to be the first branching time  and $y:=X_1(\tau)$ the first branching location. Define the events
	\begin{equation}
		T_{A}:=\left\{X\,|\,\tau\in A \right\},\quad
		L_{B}:=\left\{X\,|\, y\in B   \right\},
	\end{equation}
	where $A\subset[0,t]$ and $B\subset(-\infty,\infty)$. We estimate probabilities of the form \footnote{Note that although we consider this joint event in \eqref{prob_general_1}, asymptotics for the conditional probability 
$
		\Prob{X\in T_{A}\cap L_{B}\,|\,X_{max}(t)\leq \sqrt{2}\alpha t}
$	can be obtained from our results, since the asymptotic behavior of $\Prob{X_{max}(t)\leq \sqrt{2}\alpha t}$ is known  \cite{derrida2016large}.}

	\begin{equation}\label{prob_general_1}
	\Prob{X_{max}(t)\leq \sqrt{2}\alpha t,\,X\in T_{A}\cap L_{B}}.
	\end{equation}

	 First, we give a large deviation estimate for the probability that the maximum of a BBM is below $\sqrt{2}\a t$, with  $\alpha\in(-\infty,1)$, $\tau\in [0,\gamma t]$, and $\gamma\in(0,1]$.
	
	\begin{theorem}\label{time}
		For all $\alpha\in\left(-\infty,1\right)$ and $\gamma\in\left(0,1\right]$, 
		\begin{equation}\label{time_eq}
			\lim_{t\rightarrow\infty}\frac{1}{t}\log\Prob{X_{max}(t)\leq \sqrt{2}\alpha t,\,\,X\in T_{[0,\gamma t]}}= - \psi_1^{\alpha}\left(\gamma\right),
		\end{equation}
		where
		\begin{equation}\label{psi1}
			\psi_1^{\alpha}\left(\gamma\right):=\begin{cases}
			-\gamma+\frac{2\alpha^2}{1+\gamma}+2,&\text{ if }\gamma\in\left(0,-\frac{\alpha+\rho}{\rho}\wedge 1  \right],\\
			-(4\sqrt{2}\rho-1)\gamma+4\rho(1-\alpha),&\text{ if  }\gamma\in\left(-\frac{\alpha+\rho}{\rho}\wedge 1,\frac{1-\alpha}{2\sqrt{2}-1}\wedge 1\right],\\
			\gamma+\frac{(\alpha-(1-\gamma))^2}{\gamma},&\text{ if  }\gamma\in\left(\frac{1-\alpha}{2\sqrt{2}-1}\wedge 1,\frac{1-\alpha}{\sqrt{2}}\wedge 1\right],\\
			2\rho(1-\alpha),&\text{ if }\gamma\in\left(\frac{1-\alpha}{\sqrt{2}}\wedge 1,1\right].
			\end{cases} 
		\end{equation}
	\end{theorem}

	\noindent Note that not all cases in \eqref{time_eq} occur for  $\alpha$ fixed, as
	\begin{align}\label{int.2}
	-\frac{\alpha+\rho}{\rho}> 0 &\iff \alpha< -\rho,\quad
	-\frac{\alpha+\rho}{\rho}< 1\iff \alpha> -2\rho,\\
	\frac{1-\alpha}{2\sqrt{2}-1}<1&\iff \alpha>-2\rho,\quad
	\frac{1-\alpha}{\sqrt{2}}< 1\iff \alpha>-\rho.
	\end{align}
	
\begin{figure}[!h]
		\begin{subfigure}{.45\textwidth}
			\centering
			\includegraphics[width=.81\linewidth]{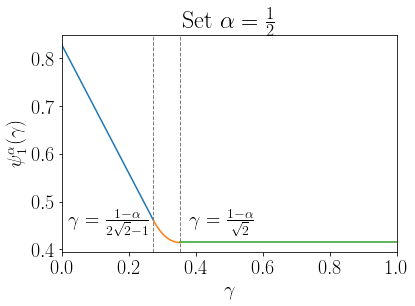}  
		\end{subfigure}
		\begin{subfigure}{.45\textwidth}
			\centering
			\includegraphics[width=.81\linewidth]{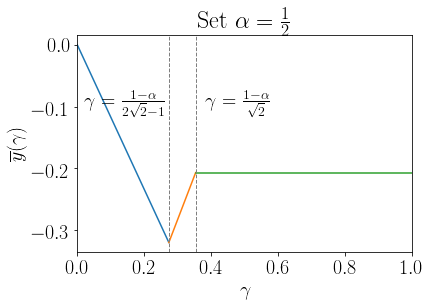}  
		\end{subfigure}
	
		\vspace{.3cm}
		\begin{subfigure}{.45\textwidth}
			\centering
			\includegraphics[width=.81\linewidth]{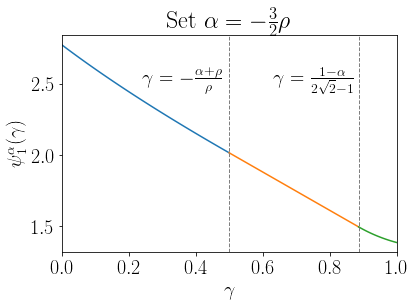}  
		\end{subfigure}
		\begin{subfigure}{.45\textwidth}
			\centering
			\includegraphics[width=.81\linewidth]{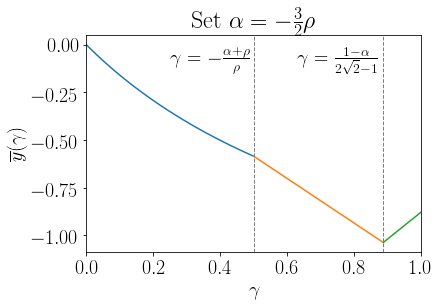}  
		\end{subfigure}
	
		\vspace{.3cm}
		\begin{subfigure}{.45\textwidth}
			\centering
			\includegraphics[width=.81\linewidth]{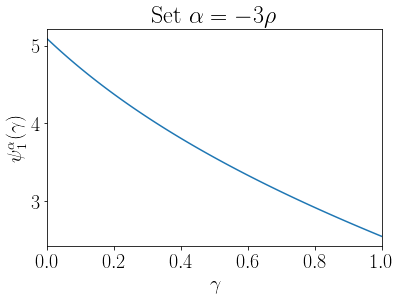}  
		\end{subfigure}
		\begin{subfigure}{.45\textwidth}
			\centering
			\includegraphics[width=.81\linewidth]{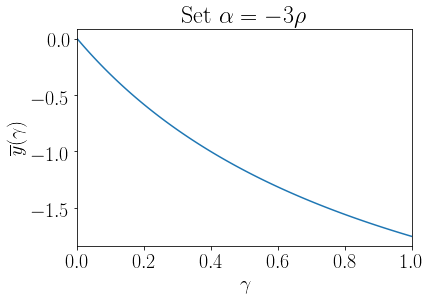}  
		\end{subfigure}
		\caption{The plots in the left column show the rate function $\psi_1^{\alpha}(\gamma)$ in \eqref{psi1} in Theorem \ref{time}. Correspondingly, the plots in the right column depict the normalized optimal first branching locations $\overline{y}(\gamma):=\lim_{t\rightarrow\infty}y(\gamma)/t$, where $y(\gamma)$ is recorded in \eqref{time_y}. We choose one representative $\alpha$ value in each of the ranges $[-\rho,1)$, $(-2\rho,-\rho)$, and $(-\infty,-2\rho]$.
		}
		\label{time_plot}
	\end{figure}

	\noindent The proof of Theorem \ref{time} in Section \ref{proof_time} shows that the optimal strategy is to make the first branching happen at time
	\begin{equation}\label{time_tau}
	\tau(\gamma)= \left(\gamma \wedge \frac{1-\alpha}{\sqrt{2}}\right) t+o(t),
	\end{equation}
	and at position
	\begin{equation}\label{time_y}
	y(\gamma)=\begin{cases}
	\frac{2\sqrt{2}\alpha \gamma}{1+\gamma}t+o(t),&\text{ if }\gamma\in\left(0,-\frac{\alpha+\rho}{\rho}\wedge 1 \right],\\
	-2\sqrt{2}\rho\gamma t+o(t),&\text{ if  }\gamma\in\left(-\frac{\alpha+\rho}{\rho}\wedge 1,\frac{1-\alpha}{2\sqrt{2}-1}\wedge 1\right],\\
	\sqrt{2}\alpha t-\sqrt{2}(1-\gamma)t+o(t),&\text{ if  }\gamma\in\left(\frac{1-\alpha}{2\sqrt{2}-1}\wedge 1,\frac{1-\alpha}{\sqrt{2}}\wedge 1\right],\\
	-\rho(1-\alpha)t+o(t),&\text{ if }\gamma\in\left(\frac{1-\alpha}{\sqrt{2}}\wedge 1,1\right].
	\end{cases}
	\end{equation}
	See Figure \ref{time_plot} for plots of $\psi_1^{\alpha}(\gamma)$ and $y(\gamma)$ as illustrations. By comparing the time-constrained optimal choices of $\tau$ in \eqref{time_tau} with the unconstrained ones in \eqref{ds_optimal_tau}, we see that to obtain a low maximum the first branching happens as late as possible, until the unrestricted optimal branching time is smaller than $\gamma t$. The time-constrained optimal choices for $y$, as shown in \eqref{time_y}, depend  on the values of $\alpha$ and $\gamma$.

	Next, we impose the restriction $X\in T_{[\gamma t, t]}$ with $\gamma\in\left(\frac{1-\alpha}{\sqrt{2}},1\right]$ to understand how an unusually late first branching time affects the large deviation estimates when $\alpha\in(-\rho,1)$. This is done in the following theorem.
	
	\begin{theorem}\label{time_late}
		For all $\alpha\in(-\rho,1)$ and $\gamma\in\left(\frac{1-\alpha}{\sqrt{2}},1\right]$,
		\begin{equation}\label{time_late_thm_eq}
			\lim_{t\rightarrow\infty}\frac{1}{t}\log\Prob{X_{max}(t)\leq \sqrt{2}\alpha t,\,\,X\in T_{[\gamma t,t]}}= -\psi_2^{\alpha}\left(\gamma\right),
		\end{equation}
		where
		\begin{equation}\label{psi2}
		\psi_2^{\alpha}\left(\gamma\right):=
			\begin{cases}
			\gamma+\frac{(\alpha-(1-\gamma))^2}{\gamma},&\text{ if  }\gamma\in\left(\frac{1-\alpha}{\sqrt{2}},(1-\alpha)\wedge 1\right),\\
			\gamma,&\text{ if }\gamma\in\left[(1-\alpha)\wedge 1,1\right].
			\end{cases}
		\end{equation}
	\end{theorem}	

\begin{figure}[!h]
		\begin{subfigure}{.45\textwidth}
			\centering
			\includegraphics[width=.81\linewidth]{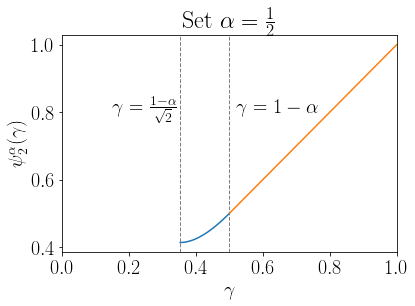}  
		\end{subfigure}
		\begin{subfigure}{.45\textwidth}
			\centering
			\includegraphics[width=.81\linewidth]{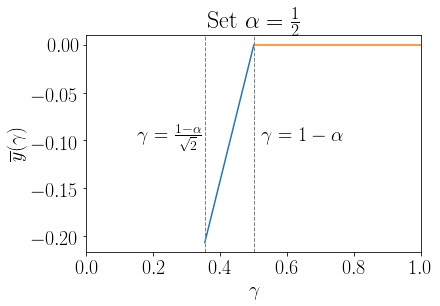}  
		\end{subfigure}
		
		\vspace{.3cm}
		\begin{subfigure}{.45\textwidth}
			\centering
			\includegraphics[width=.81\linewidth]{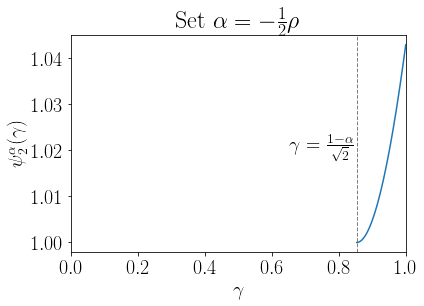}  
		\end{subfigure}
		\begin{subfigure}{.45\textwidth}
			\centering
			\includegraphics[width=.81\linewidth]{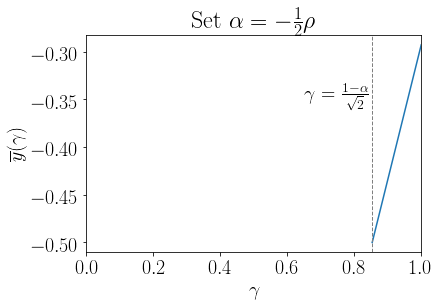}  
		\end{subfigure}

		\caption{The plots in the left column show the rate function $\psi_2^{\alpha}(\gamma)$ in \eqref{psi2} in Theorem \ref{time_late}. Correspondingly, the plots in the right column depict the normalized optimal first branching locations $\overline{y}(\gamma)=\lim_{t\rightarrow\infty}y(\gamma)/t$, where $y(\gamma)$ is recorded in \eqref{time_late_y}. We choose one representative $\alpha$ value in each of the two ranges $(0,1)$ and $(-\rho,0]$. 
		}
		\label{time_late_plot}
	\end{figure}

	\noindent The proof of Theorem \ref{time_late} in Section \ref{proof_time} shows that the optimal strategy is to let the first branching happens at time
	\begin{equation}\label{time_late_tau}
		\tau(\gamma)=\gamma t+o(t), \text{ for all }\gamma\in\left(\frac{1-\alpha}{\sqrt{2}},1\right],
	\end{equation}
	and at position
	\begin{equation}\label{time_late_y}
		y(\gamma)=\begin{cases}
		\sqrt{2}\alpha t-\sqrt{2}(1-\gamma)t+o(t),&\text{ if  }\gamma\in\left(\frac{1-\alpha}{\sqrt{2}},(1-\alpha)\wedge 1\right),\\
		o(t),&\text{ if }\gamma\in\left[(1-\alpha)\wedge 1,1\right].
		\end{cases}
	\end{equation}
	See Figure \ref{time_late_plot} for plots of $\psi_{2}^{\alpha}(\gamma)$ and $y(\gamma)$ as illustrations. Note that if $\gamma>1-\alpha$, the probability with the late-first-branching constraint is of order $\eee^{-\gamma t+o(t)}$, which is of the same order as the probability that a  BBM does not  branch in $[0,\g t]$. 
	
	After studying how restrictions on the first branching time affect the large deviation estimates, we turn our attention to the effects of a constrained first branching location. In Theorems \ref{location_below} and \ref{location_above}, we fix the first branching time to be in the interval $[(\gamma-\epsilon)t,\gamma t]$ where $\epsilon$ is positive and small, and impose restrictions on the first branching location to be either below or above $\sqrt{2}\alpha t-\sqrt{2}(t-\tau)$.
	
	 Note that $\sqrt{2}(t-\tau)$ is the leading order of the maximum of a BBM running for time $t-\tau$. If the maximum of the BBM  (running for time $t$) has to stay below $\sqrt{2}\alpha t$ and 
	  the first branching location is below  $\sqrt{2}\alpha t-\sqrt{2}(t-\tau)$, the two BBMs starting from the initial branching position do not need to have an unusually low maxima. 
	  The following theorem gives the large deviation estimates when the first branching location is forced to be below $\sqrt{2}\alpha t-\sqrt{2}(t-\tau)$. 
	
	\begin{theorem}\label{location_below}
		For all $\alpha\in(-\infty,1)$, $\gamma\in\left(0,1\right]$, and $\beta\in[1,\infty)$, 
		\begin{equation}\label{location_below_eq}
			\lim_{\epsilon\rightarrow 0}\lim_{t\rightarrow\infty}\frac{1}{t}\log \Prob{X_{max}(t)\leq \sqrt{2}\alpha t,\,X\in T_{[(\gamma-\epsilon)t,\gamma t]}\cap L_{\left(-\infty,\sqrt{2}\alpha t-\sqrt{2}\beta(t-\tau)\right]}}=-\psi_3^{\alpha,\gamma}\left(\beta\right),
		\end{equation}
		where
		\begin{equation}\label{psi3}
		\psi_3^{\alpha,\gamma}\left(\beta\right):=
		\begin{cases}
		\gamma,&\text{ if }\beta\in\left[1,\beta_1^{\alpha}(\gamma)\vee 1\right],\\
		\gamma+\frac{(\alpha-(1-\gamma)\beta)^2}{\gamma},&\text{ if }\beta\in\left(\beta_1^{\alpha}(\gamma)\vee 1,\infty\right),
		\end{cases}  \text{ with }
		\beta_1^{\alpha}(\gamma):=\begin{cases}
		\frac{\alpha}{1-\gamma}, &\text{ if }0<\gamma<1,\\
		+\infty,&\text{ if }\gamma=1.
		\end{cases}
		\end{equation}
	\end{theorem}
	
\begin{figure}[!h]
		\begin{subfigure}{.45\textwidth}
			\centering
			\includegraphics[width=.81\linewidth]{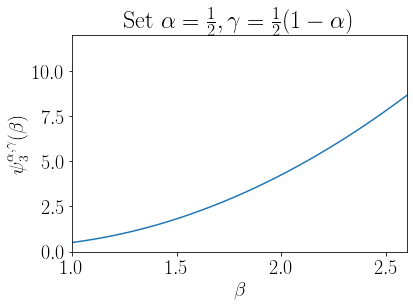}  
		\end{subfigure}
		\begin{subfigure}{.45\textwidth}
			\centering
			\includegraphics[width=.81\linewidth]{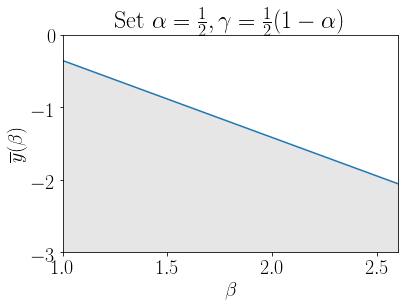}  
		\end{subfigure}
		
		\vspace{.3cm}
		\begin{subfigure}{.45\textwidth}
			\centering
			\includegraphics[width=.81\linewidth]{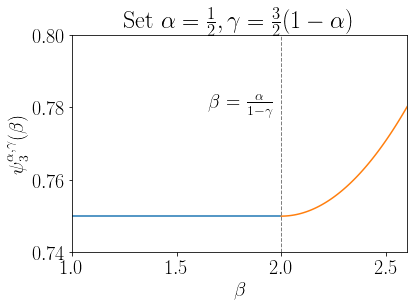}  
		\end{subfigure}
		\begin{subfigure}{.45\textwidth}
			\centering
			\includegraphics[width=.81\linewidth]{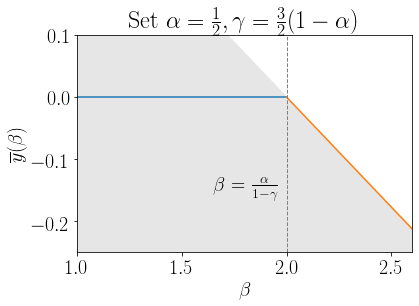}  
		\end{subfigure}
		
		\caption{The plots in the left column show the rate function $\psi_3^{\alpha,\gamma}(\beta)$ in \eqref{psi3} in Theorem \ref{location_below}. Correspondingly, the plots in the right column depict the normalized optimal first branching locations $\overline{y}(\beta)=\lim_{t\rightarrow\infty}y(\beta)/t$, where $y(\beta)$ is recorded in \eqref{location_below_y}. Representative values are chosen for $\alpha$ and $\gamma$. The shaded areas in the right plots indicate the area where the first branching is allowed to happen.
		}
		\label{location_below_plot}
	\end{figure}

	\noindent The proof of Theorem \ref{location_below} in Section \ref{proof_location} shows that the optimal strategy for the initial particle is to split at location
	\begin{equation}\label{location_below_y}
		y(\beta)=\begin{cases}
		o(t),&\text{ if }\beta\in\left[1,\beta_1^{\alpha}(\gamma)\vee 1\right],\\
		\sqrt{2}\alpha t-\sqrt{2}\beta(1-\gamma)t+o(t),&\text{ if }\beta\in\left(\beta_1^{\alpha}(\gamma)\vee 1,\infty\right).
		\end{cases}
	\end{equation}
	See Figure \ref{location_below_plot} for plots of $\psi_{3}^{\alpha,\gamma}(\beta)$ and $y(\beta)$ as illustrations. Note that since $\beta_1^{\alpha}(\gamma)>1\iff \gamma>1-\alpha$, the $\beta\in\left[1,\beta_1^{\alpha}(\gamma)\vee 1\right]$ case occurs only when $\alpha\in(0,1)$ and $\gamma\in(1-\alpha,1]$. In this situation, the additional location restriction has no additional impact and the values of the rate functions $\psi_{2}^{\alpha}(\gamma)$ and $\psi_{3}^{\alpha,\gamma}(\beta)$ agree.
	In all other cases, the location constraint does affect the large deviation estimates, and the optimal strategy for the initial particle is to split at the highest possible position $y\sim\sqrt{2}t-\sqrt{2}\beta(1-\gamma)t$.  

	Next, we consider the case where the first branching location is restricted to be above $\sqrt{2}\alpha t-\sqrt{2}(t-\tau)$. 	
	\begin{theorem}\label{location_above}
		For all $\alpha\in\left(-\infty,1\right)$, $\gamma\in\left(0,1\right]$, and $\beta\in(-\infty,1]$, 
		\begin{equation}\label{location_above_eq}
			\lim_{\epsilon\rightarrow 0}\lim_{t\rightarrow\infty}\frac{1}{t}\log \Prob{X_{max}(t)\leq \sqrt{2}\alpha t,\,X\in T_{[(\gamma-\epsilon)t,\gamma t]}\cap L_{\left[\sqrt{2}\alpha t-\sqrt{2}\beta(t-\tau),\infty\right)}}=-\psi_4^{\alpha,\gamma}\left(\beta \right),
		\end{equation}
		where,  
				\begin{itemize}
			\item  for $\gamma\in\left(0,-\frac{\alpha+\rho}{\rho}\wedge 1 \right)$, $\psi_4^{\alpha,\gamma}\left(\beta \right)$ is equal to
			\begin{equation}\label{location_above_1}
			 \begin{split}
			\begin{cases}
			-\left(1+\beta^2 \right)\gamma+\frac{(\alpha-\beta)^2}{\gamma}+2(\alpha\beta+1),&\text{ if }\beta\in\left(-\infty,\frac{\alpha}{1+\gamma}\right),\\
			-\gamma+\frac{2\alpha^2}{1+\gamma}+2,&\text{ if }\beta\in\left[\frac{\alpha}{1+\gamma},1 \right];
			\end{cases} 			
			\end{split}
			\end{equation}
			
			\item for $\gamma\in\left[-\frac{\alpha+\rho}{\rho}\wedge 1, 1\right]$, $\psi_4^{\alpha,\gamma}\left(\beta \right)$ is equal to
			\begin{equation}\label{location_above_2}
			\begin{cases}
			-\left(1+\beta^2 \right)\gamma+\frac{(\alpha-\beta)^2}{\gamma}+2(\alpha\beta+1) ,&\text{ if }\beta\in(-\infty,-\rho],\\
			-\left(4\rho(1-\beta)-1-\beta^2 \right)\gamma+\frac{(\alpha-\beta)^2}{\gamma}+4\rho(1-\beta)+2(\alpha-\beta)\beta,&\text{ if }\beta\in\left(-\rho, \beta_2^{\alpha}(\gamma)\wedge 1 \right),\\
			-\left(4\sqrt{2}\rho-1 \right)\gamma+4\rho(1-\alpha),&\text{ if }\beta\in\left[\beta_2^{\alpha}(\gamma)\wedge 1,1 \right],			
			\end{cases}
			\end{equation}
		\end{itemize}
		with
		\begin{equation}
			\beta_2^{\alpha}(\gamma):=\begin{cases}
			\frac{\alpha+2\rho\gamma}{1-\gamma},&\text{ if }0<\gamma<1,\\
			-\infty,&\text{ if }\gamma=1\text{ and }\alpha<-2\rho,\\
			+\infty,&\text{ if }\gamma=1\text{ and }\alpha\geq-2\rho.
			\end{cases}
		\end{equation}
	\end{theorem}

	\begin{figure}[!t]
		\begin{subfigure}{.45\textwidth}
			\centering
			\includegraphics[width=.81\linewidth]{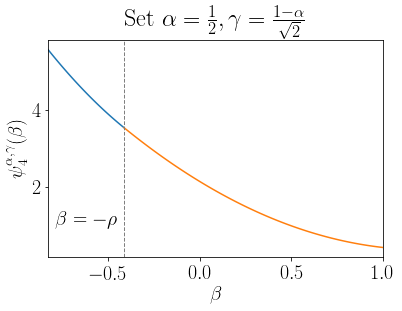}  
		\end{subfigure}
		\begin{subfigure}{.45\textwidth}
			\centering
			\includegraphics[width=.81\linewidth]{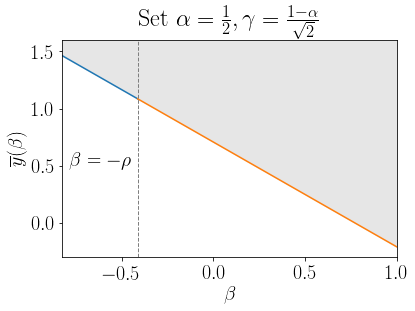}  
		\end{subfigure}
		
		\vspace{.3cm}
		\begin{subfigure}{.45\textwidth}
			\centering
			\includegraphics[width=.81\linewidth]{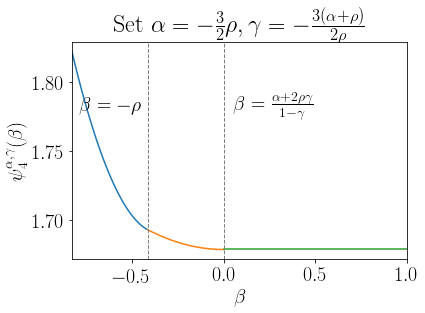}  
		\end{subfigure}
		\begin{subfigure}{.45\textwidth}
			\centering
			\includegraphics[width=.81\linewidth]{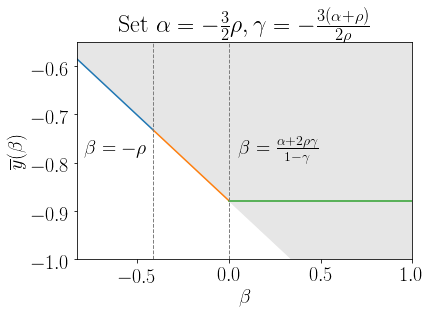}  
		\end{subfigure}
		
		\vspace{.3cm}
		\begin{subfigure}{.45\textwidth}
			\centering
			\includegraphics[width=.81\linewidth]{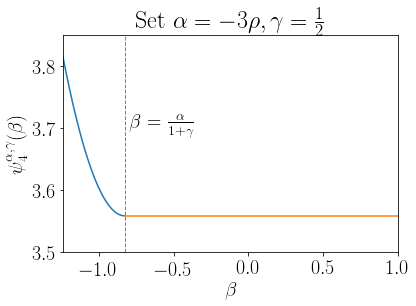}  
		\end{subfigure}
		\begin{subfigure}{.45\textwidth}
			\centering
			\includegraphics[width=.81\linewidth]{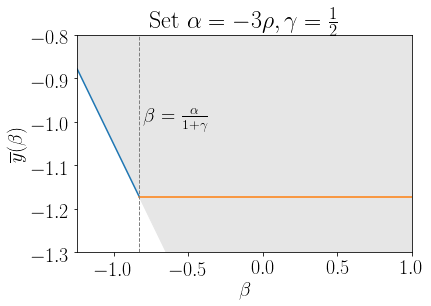}  
		\end{subfigure}
		
		\caption{The plots in the left column show the rate function $\psi_4^{\alpha,\gamma}(\beta)$ in \eqref{location_above_1} and \eqref{location_above_2} in Theorem \ref{location_below}. Correspondingly, the plots in the right column depict the normalized optimal first branching locations $\overline{y}(\beta)=\lim_{t\rightarrow\infty}y(\beta)/t$, where $y(\beta)$ is recorded in \eqref{location_above_1_y} and \eqref{location_above_2_y}. Representative values are chosen for $\alpha$ and $\gamma$. The shaded areas in the right plots indicate again the area where the first branching is allowed to happen.}
		\label{location_above_plot}
	\end{figure}

	\noindent The proof of Theorem \ref{location_above} in Section \ref{proof_location} shows that the optimal first branching location for $\gamma\in\left(0,-\frac{\alpha+\rho}{\rho}\wedge 1 \right)$ is given by
		\begin{equation}\label{location_above_1_y}
		y(\beta)=\begin{cases}
		\sqrt{2}\alpha t-\sqrt{2}\beta(1-\gamma) t+o(t),&\text{ if }\beta\in\left(-\infty,\frac{\alpha}{1+\gamma}\right),\\
		\frac{2\sqrt{2}\alpha\gamma }{1+\gamma}t+o(t),&\text{ if }\beta\in\left[\frac{\alpha}{1+\gamma},1 \right],
		\end{cases}
		\end{equation}
	while for  $\gamma\in\left[-\frac{\alpha+\rho}{\rho}\wedge 1, 1\right]$, it is equal to
		\begin{equation}\label{location_above_2_y}
		y(\beta)=\begin{cases}
		\sqrt{2}\alpha t-\sqrt{2}\beta(1-\gamma) t+o(t),&\text{ if }\beta\in\left(-\infty,\beta_2^{\alpha}(\gamma)\wedge 1 \right),\\
		-2\sqrt{2}\rho\gamma t+o(t),&\text{ if }\beta\in\left[\beta_2^{\alpha}(\gamma)\wedge 1,1 \right].			
		\end{cases}
		\end{equation}
	See Figure \ref{location_above_plot} for plots of $\psi_{4}^{\alpha,\gamma}(\beta)$ and $y(\beta)$ as illustrations. Note that, by \eqref{int.2}, for $\alpha\in(-\infty,-2\rho]$ and for $\alpha\in[-\rho,1)$, \eqref{location_above_1} resp. \eqref{location_above_2} define the rate function for all $\gamma\in(0,1)$, while for $\alpha\in(-2\rho,-\rho)$, both \eqref{location_above_1} and \eqref{location_above_2} apply in appropriate $\gamma$ ranges. Moreover, as $\beta_2^{\alpha}(\gamma)=\frac{\alpha+2\rho\gamma}{1-\gamma}$ when $0<\gamma<1$ and
	\begin{equation}
	\frac{\alpha+2\rho\gamma}{1-\gamma}< 1\iff \gamma<\frac{1-\alpha}{2\sqrt{2}-1},
	\end{equation}
	the last case of \eqref{location_above_2} only occurs when the first branching happens early enough. It is also worth noticing that, if the range of $\beta$ includes 1 (see the second case in \eqref{location_above_1} and the last two cases in \eqref{location_above_2}), the location constraint in Theorem \ref{location_above} has no additional impact and the rate functions $\psi_{4}^{\alpha,\gamma}(\beta)$ and $\psi_{1}^{\alpha,\gamma}(\beta)$ coincide.
	 In all other cases, the effect of the location constraint is evident, and the optimal strategy requires the initial particle to split at the lowest possible position $y\sim \sqrt{2}\alpha t-\sqrt{2}\beta(1-\gamma)t$ to minimize such an effect.
	 
	 \paragraph{Outline of the paper. } In Section 2, we  first use the branching property to decompose the probability of a BBM to have a maximum below $\sqrt{2}\alpha t$ with constrained first branching time and location. The decomposition is stated in  Lemma \ref{decomposition}. In Lemmas \ref{y1} to \ref{y3}, we analyze the resulting terms of this decomposition separately. In Section 3 we then prove our main results Theorems \ref{time} to \ref{location_above}, based on the preparatory lemmas from Section 2.

%


\section{Preparatory Lemmas}\label{preliminary}

	
	\subsection{Decomposition of the refined large deviation probabilities}
	
	We rewrite the estimates from \cite{derrida2016large} (see also \eqref{ds}) in a form that is convenient for us.
	
	\begin{corollary}\label{prob}
		For any $\tau\in\left[0,t\right)$, $\alpha,y\in\R$, as $t\rightarrow\infty$, 
		\begin{numcases}{\Prob{X_{max}(t-\tau)\leq\sqrt{2}\alpha t-y}=}
		1-\eee^{-O(t)}, & if $y\in I_1$,\label{prob_1}\\
		\mathrm{e}^{-\sqrt{2}\rho
			\left(\sqrt{2}\left(1-\alpha\right)t-\sqrt{2}\tau+y\right)+o(t)},& if $y\in I_2$,\label{prob_2}\\
		\mathrm{e}^{-(t-\tau)-\frac{(\sqrt{2}\alpha t-y)^2}{2(t-\tau)}+o(t)},& if $y\in I_3$,\label{prob_3}
		\end{numcases}
		where
		\begin{equation}
		\begin{split}
		I_1&:=\left(-\infty, \sqrt{2}\alpha t-\sqrt{2}\left(t-\tau\right)\right),\\
		I_2&:=\left[\sqrt{2}\alpha t-\sqrt{2}\left(t-\tau\right),\sqrt{2}\alpha t+\sqrt{2}\rho\left(t-\tau\right)\right],\\
		I_3&:=\left( \sqrt{2}\alpha t+\sqrt{2}\rho\left(t-\tau\right),\infty\right).
		\end{split}
		\end{equation}
	\end{corollary}
	
	\begin{proof}
		By rewriting the probability in the corollary as 
		\begin{equation}\label{prob_rewrite}
		\Prob{X_{max}(t-\tau)\leq\sqrt{2}\,\,\frac{  \sqrt{2}\alpha t-y}{\sqrt{2}(t-\tau)}\,\,(t-\tau)},
		\end{equation}
		applying the estimate \eqref{ds}, and distinguishing whether $\frac{  \sqrt{2}\alpha t-y}{\sqrt{2}(t-\tau)}$ is in the range $[-\rho,1]$ or $(-\infty,-\rho]$, \eqref{prob_2} and \eqref{prob_3} follow. For \eqref{prob_1}, we apply \eqref{u_3} with
		\begin{equation}
			x=\sqrt{2}\alpha t-\sqrt{2}(t-\tau)-y
		\end{equation}
	to  one minus the probability in \eqref{prob_rewrite}, obtaining that one minus the probability in \eqref{prob_rewrite} is asymptotically equal to
		\begin{equation}
			\begin{split}
			&C\left(t-\tau \right)^{-\frac{1}{2}}\exp(-\sqrt{2}\left(\sqrt{2}\alpha t-\sqrt{2}(t-\tau)-y\right)-\frac{\left(\sqrt{2}\alpha t-\sqrt{2}(t-\tau)-y\right)^2}{2(t-\tau)}),
			\end{split}
		\end{equation}
		which is equal to $\eee^{-O(t)}$.  This implies \eqref{prob_1}.
	\end{proof}
	
	Recall that in Theorems \ref{time} to \ref{location_above}, the refined large deviation probabilities take the form
	\begin{equation}\label{prob_general}
	\Prob{X_{max}(t)\leq \sqrt{2}\alpha t,\,X\in T_{A}\cap L_{B}},
	\end{equation}
	where $A\subset[0,t]$ and $B\subset(-\infty,\infty)$. We rewrite \eqref{prob_general} by disintegrating at the first branching time and using the branching property of BBM, which results in the following lemma.

	
	\begin{lemma}\label{decomposition}
		For all $\alpha\in(-\infty,1)$, $A\subset[0,t]$, and $B\subset (-\infty,\infty)$, as $t\rightarrow\infty$,
		\begin{equation}\label{decomposition_eq}
			\Prob{X_{max}(t)\leq \sqrt{2}\alpha t,\,X\in T_{A}\cap L_{B}}=\int_{A}\left(Y_1(B\cap I_1)+Y_2(B\cap I_2)+Y_3(B\cap I_3)\right)\dd\tau,
		\end{equation}
		where
		\begin{align}
			Y_1(B\cap I_1)&:=\int_{B\cap I_1}\mathrm{e}^{-\tau}\frac{1}{\sqrt{2\pi\tau}}\exp(-\frac{y^2}{2\tau})\dd y,\\
			Y_2(B\cap I_2)&:=\int_{B\cap I_2}\mathrm{e}^{\left(4\sqrt{2}\rho-1\right)\tau-4\rho(1-\alpha)t+o(t)}\frac{1}{\sqrt{2\pi\tau}}\exp(-\frac{(y+2\sqrt{2}\rho\tau)^2}{2\tau})\dd y,\label{eq.y2}\\
			Y_3(B\cap I_3)&:=\int_{B\cap I_3}\sqrt{\frac{t-\tau}{t+\tau}}\mathrm{e}^{-2t+\tau-\frac{2\alpha^2t^2}{t+\tau}+o(t)}\frac{1}{\sqrt{2\pi\tau\frac{t-\tau}{t+\tau}}}\exp(-\frac{(y-\frac{2\sqrt{2}\alpha t\tau}{t+\tau})^2}{2\tau\frac{t-\tau}{t+\tau}})\dd y.
		\end{align}
	\end{lemma}
	
	\begin{proof}
		By  the branching property of BBM at the first branching time $\tau$ and location $y$, the probability in \eqref{decomposition_eq} is equal to
		\begin{equation}
			\int_{A}\int_{B}\frac{1}{\sqrt{2\pi\tau}}\mathrm{e}^{-\tau-\frac{y^2}{2\tau}}\Prob{X_{max}(t-\tau)\leq\sqrt{2}\alpha t-y}^2\dd y\dd\tau,
		\end{equation}
		which by Corollary \ref{prob} is
		\begin{equation}
			\begin{split}
			&\int_{A}\int_{B\cap I_1}\left(1-\eee^{-O(t)}\right)\frac{1}{\sqrt{2\pi\tau}}\mathrm{e}^{-\tau-\frac{y^2}{2\tau}}\dd y\dd\tau\\
			+\,\,&\int_{A}\int_{B\cap I_2}\frac{1}{\sqrt{2\pi\tau}}\mathrm{e}^{-\tau-\frac{y^2}{2\tau}-2\sqrt{2}\rho
				\left(\sqrt{2}\left(1-\alpha\right)t-\sqrt{2}\tau+y\right)+o(t)}\dd y\dd\tau\\
			+\,\,&\int_{A}\int_{B\cap I_3}\frac{1}{\sqrt{2\pi\tau}}\mathrm{e}^{-\tau-\frac{y^2}{2\tau}-2(t-\tau)-\frac{(\sqrt{2}\alpha t-y)^2}{(t-\tau)}+o(t)}\dd y\dd\tau.
			\end{split}
		\end{equation}
		Completing the squares for $y$, rearranging the terms, and letting $t\rightarrow\infty$, we obtain \eqref{decomposition_eq}.
	\end{proof}
\subsection{First estimates on $Y_1$, $Y_2$ and $Y_3$}
	The decomposition in Lemma \ref{decomposition} suggests 
	we need to obtain good approximations for  $Y_1(B\cap I_1), Y_2(B\cap I_2), $ and $Y_3(B\cap I_3)$, which is done in Lemma \ref{y1}. Lemma \ref{y2}, and Lemma \ref{y3}.


%
	
	\begin{lemma}\label{y1}
		For all $t>0$ and $\alpha\in(-\infty,1)$, suppose there exist some $\lambda_\tau\in(0,1]$ and $\lambda_u\in\R$ such that $\tau=\lambda_\tau t$ and $u=\sqrt{2}\alpha t-\sqrt{2}\lambda_u(1-\lambda_\tau)t$. Then, as $t\rightarrow\infty$,
		\begin{equation}\label{y1_eq}
			Y_1((-\infty,u])=\begin{cases}
			\exp(-t\ll{11}{\lambda_\tau}{\lambda_u}+o(t)),&\text{ if }\lambda_u>\beta_1^{\alpha}(\lambda_\tau),\\
			\exp(-t\l{12}{\lambda_\tau}+o(t)), &\text{ if } \lambda_u\leq\beta_1^{\alpha}(\lambda_\tau),
			\end{cases}			
		\end{equation}
		where
		\begin{equation}
			\begin{split}
			\ll{11}{x}{y}:=x+\frac{\left(\alpha-y(1-x) \right)^2}{x},\quad
			\l{12}{x}:=x.
			\end{split}
		\end{equation}
	\end{lemma}

	\begin{proof}
		Since $\tau=\lambda_\tau t$ and $u=\sqrt{2}\alpha t-\sqrt{2}\lambda_u(1-\lambda_\tau)t$, $Y_1((-\infty,u])$ is equal to
		\begin{equation}\label{y1_eq_pf}
			\int_{-\infty}^{\sqrt{2}\alpha t-\sqrt{2}\lambda_u(1-\lambda_\tau)t}\mathrm{e}^{-\lambda_\tau t}\frac{1}{\sqrt{2\pi\lambda_\tau t}}\exp(-\frac{y^2}{2\lambda_\tau t})\dd y.
		\end{equation}
		When $\sqrt{2}\alpha t-\sqrt{2}\lambda_u(1-\lambda_\tau)t<0\iff \lambda_u\left(1-\lambda_\tau\right)>\alpha$, by Gaussian tail estimates \eqref{y1_eq_pf} is approximately
		\begin{equation}
			\exp(-\lambda_\tau t-\frac{(\alpha-\lambda_u(1-\lambda_\tau))^2}{\lambda_\tau}t+o(t)).
		\end{equation} 
		On the other hand, when $\lambda_u\left(1-\lambda_\tau\right)\leq\alpha$, the $y$ integral in \eqref{y1_eq_pf} is bounded from below by $\frac{1}{2}$ and from above by $1$. Combining these two cases, \eqref{y1_eq} follows.
	\end{proof}

	\begin{lemma}\label{y2}
		For all $t>0$ and $\alpha\in(-\infty,1)$, suppose there exist some $0<\lambda_\tau\leq 1$ and $-\infty<\lambda_v<\lambda_u<\infty$, such that $\tau=\lambda_\tau t$, $u=\sqrt{2}\alpha t-\sqrt{2}\lambda_u(1-\lambda_\tau)t$, and $v=\sqrt{2}\alpha t-\sqrt{2}\lambda_v(1-\lambda_\tau)t$. Then, as $t\rightarrow\infty$,
		\begin{equation}\label{y2_eq}
		Y_2([u,v])=\begin{cases}
		\exp(-t\ll{21}{\lambda_\tau}{\lambda_v}+o(t)),&\text{ if }\lambda_v>\beta_2^{\alpha}(\lambda_\tau),\\
		\exp(-t\l{22}{\lambda_\tau}+o(t)), &\text{ if }\lambda_v\leq\beta_2^{\alpha}(\lambda_\tau)\leq\lambda_u,\\
		\exp(-t\ll{21}{\lambda_\tau}{\lambda_u}+o(t)),&\text{ if }\lambda_u<\beta_2^{\alpha}(\lambda_\tau),
		\end{cases}
		\end{equation}
		where
		\begin{equation}
			\begin{split}
			\ll{21}{x}{y}&:=-\left(4\rho(1-y)-1-y^2\right)x+\frac{(\alpha-y)^2}{x}+\left(4\rho\left(1-y\right)+2y\left(\alpha-y\right)\right),\\
			\l{22}{x}&:=-\left(4\sqrt{2}\rho-1 \right)x+4\rho\left(1-\alpha\right).
			\end{split}
		\end{equation}
	\end{lemma}

	\begin{proof}
		Plugging  $\tau=\lambda_\tau t$, $u=\sqrt{2}\alpha t-\sqrt{2}\lambda_u(1-\lambda_\tau)t$, and $v=\sqrt{2}\alpha t-\sqrt{2}\lambda_v(1-\lambda_\tau)t$ into \eqref{eq.y2}, $Y_2([u,v])$ becomes
		\begin{equation}
			\int_{\sqrt{2}\alpha t-\sqrt{2}\lambda_u(1-\lambda_\tau)t}^{\sqrt{2}\alpha t-\sqrt{2}\lambda_v(1-\lambda_\tau)t}\mathrm{e}^{\left(4\sqrt{2}\rho-1\right)\lambda_\tau t-4\rho(1-\alpha)t+o(t)}\frac{1}{\sqrt{2\pi\lambda_\tau t}}\exp(-\frac{(y+2\sqrt{2}\rho\lambda_\tau t)^2}{2\lambda_\tau t})\dd y.
		\end{equation}
		With a change of variable $z=y+2\sqrt{2}\rho\lambda_\tau t$, this is equal to
		\begin{equation}\label{y2_eq_pf}
			\int_{\left(\sqrt{2}\alpha -\sqrt{2}\lambda_u(1-\lambda_\tau)+2\sqrt{2}\rho\lambda_\tau\right) t}^{\left(\sqrt{2}\alpha -\sqrt{2}\lambda_v(1-\lambda_\tau)+2\sqrt{2}\rho\lambda_\tau \right) t}\mathrm{e}^{\left(4\sqrt{2}\rho-1\right)\lambda_\tau t-4\rho(1-\alpha)t+o(t)}\frac{1}{\sqrt{2\pi\lambda_\tau t}}\exp(-\frac{z^2}{2\lambda_\tau t})\dd z.
		\end{equation}
		
		When $\sqrt{2}\alpha -\sqrt{2}\lambda_v(1-\lambda_\tau)+2\sqrt{2}\rho\lambda_\tau <0$,
		\eqref{y2_eq_pf} can be estimated by Gaussian tail asymptotics as
		\begin{equation}
			\exp( \left(4\sqrt{2}\rho-1\right)\lambda_\tau t-4\rho(1-\alpha)t -\frac{\left(\sqrt{2}\alpha -\sqrt{2}\lambda_v(1-\lambda_\tau)+2\sqrt{2}\rho\lambda_\tau\right)^2t}{2\lambda_\tau} +o(t)),
		\end{equation}
		which, after some rearrangements, is equal to the first term on the right-hand side of \eqref{y2_eq}.
		
		When $\sqrt{2}\alpha -\sqrt{2}\lambda_u(1-\lambda_\tau)+2\sqrt{2}\rho\lambda_\tau>0$,
		we apply  again Gaussian tail estimates to \eqref{y2_eq_pf} and obtain
		\begin{equation}
			\exp( \left(4\sqrt{2}\rho-1\right)\lambda_\tau t-4\rho(1-\alpha)t -\frac{\left(\sqrt{2}\alpha -\sqrt{2}\lambda_u(1-\lambda_\tau)+2\sqrt{2}\rho\lambda_\tau\right)^2t}{2\lambda_\tau} +o(t)),
		\end{equation} 
		which, after some rearrangements, is equal to the third term on the right-hand side of \eqref{y2_eq}.
		
		When $\sqrt{2}\alpha -\sqrt{2}\lambda_v(1-\lambda_\tau)+2\sqrt{2}\rho\lambda_\tau \geq 0$ and $\sqrt{2}\alpha -\sqrt{2}\lambda_u(1-\lambda_\tau)+2\sqrt{2}\rho\lambda_\tau  \leq 0$, 
		the $z$ integral in \eqref{y2_eq_pf} is bounded below by $\frac{1}{2}$ and above by $1$. Then \eqref{y2_eq_pf} is equal to the second term on the right-hand side of \eqref{y2_eq}.	
		
		Combining these three cases, \eqref{y2_eq} follows.	
	\end{proof}

	\begin{lemma}\label{y3}
		For all $t>0$ and $\alpha\in(-\infty,1)$, suppose there exist some $\lambda_\tau\in(0,1]$ and $\lambda_v\in\R$ such that $\tau=\lambda_\tau t$ and $v=\sqrt{2}\alpha t-\sqrt{2}\lambda_v(1-\lambda_\tau)t$. Then, as $t\rightarrow\infty$,
		\begin{equation}\label{y3_eq}
		Y_3([v,\infty))=\begin{cases}
		\exp(-t\ll{31}{\lambda_\tau}{\lambda_v}+o(t)),&\text{ if }\lambda_v<\frac{\alpha}{1+\lambda_\tau},\\
		\exp(-t\l{32}{\lambda_\tau}+o(t)), &\text{ if }\lambda_v\geq\frac{\alpha}{1+\lambda_\tau},
		\end{cases}
		\end{equation}
		where
		\begin{equation}
			\begin{split}
			\ll{31}{x}{y}&:=-\left(1+y^2 \right)x+\frac{(\alpha-y)^2}{x}+2\left(\alpha y+1\right),\\
			\l{32}{x}&:=-x+\frac{2\alpha^2}{1+x}+2.
			\end{split}
		\end{equation}
	\end{lemma}

	\begin{proof}
		With $\tau=\lambda_\tau t$ and $v=\sqrt{2}\alpha t-\sqrt{2}\lambda_v(1-\lambda_\tau)t$, $Y_3([v,\infty))$ is equal to 
		\begin{equation}
			\int_{\sqrt{2}\alpha t-\sqrt{2}\lambda_v(1-\lambda_\tau)t}^{\infty}\sqrt{\frac{1-\lambda_\tau}{1+\lambda_\tau}}\mathrm{e}^{-2t+\lambda_\tau t-\frac{2\alpha^2t}{1+\lambda_\tau }+o(t)}\frac{1}{\sqrt{2\pi\lambda_\tau t\frac{1-\lambda_\tau}{1+\lambda_\tau}}}\exp(-\frac{(y-\frac{2\sqrt{2}\alpha \lambda_\tau t}{1+\lambda_\tau})^2}{2\lambda_\tau t\frac{1-\lambda_\tau}{1+\lambda_\tau}})\dd y.
		\end{equation}
		With a change of variable $z=y-\frac{2\sqrt{2}\alpha \lambda_\tau t}{1+\lambda_\tau}$, this becomes
		\begin{equation}\label{y3_eq_pf}
			\int_{\left(\sqrt{2}\alpha -\sqrt{2}\lambda_v(1-\lambda_\tau)-\frac{2\sqrt{2}\alpha \lambda_\tau }{1+\lambda_\tau}\right)t}^{\infty}\sqrt{\frac{1-\lambda_\tau}{1+\lambda_\tau}}\mathrm{e}^{-2t+\lambda_\tau t-\frac{2\alpha^2t}{1+\lambda_\tau }+o(t)}\frac{1}{\sqrt{2\pi\lambda_\tau t\frac{1-\lambda_\tau}{1+\lambda_\tau}}}\exp(-\frac{z^2}{2\lambda_\tau t\frac{1-\lambda_\tau}{1+\lambda_\tau}})\dd z.
		\end{equation}
		
		If $\sqrt{2}\alpha -\sqrt{2}\lambda_v(1-\lambda_\tau)-\frac{2\sqrt{2}\alpha \lambda_\tau }{1+\lambda_\tau}>0$,
		by Gaussian tail estimates \eqref{y3_eq_pf} is equal to
		\begin{equation}
			\exp(-2t+\lambda_\tau t-\frac{2\alpha^2t}{1+\lambda_\tau }-\frac{\left(\sqrt{2}\alpha -\sqrt{2}\lambda_v(1-\lambda_\tau)-\frac{2\sqrt{2}\alpha \lambda_\tau }{1+\lambda_\tau}\right)^2t}{2\lambda_\tau \frac{1-\lambda_\tau}{1+\lambda_\tau}}+o(t) ),
		\end{equation}
		which, after some rearrangements, is equal to the first term on the right-hand side of \eqref{y3_eq}.
		
		For $\sqrt{2}\alpha -\sqrt{2}\lambda_v(1-\lambda_\tau)-\frac{2\sqrt{2}\alpha \lambda_\tau }{1+\lambda_\tau} \leq 0$,
		the $z$ integral in \eqref{y3_eq_pf} is bounded below by $\frac{1}{2}$ and above by $1$. After some rearrangements, \eqref{y3_eq_pf} is equal to the second term on the right-hand side of \eqref{y3_eq}.
		
		Combining these two cases, \eqref{y3_eq} follows.
	\end{proof}

	
	In addition to the estimates on $Y_1, Y_2$, and $Y_3$, to compute \eqref{decomposition_eq} in Lemma \ref{decomposition} we still need to compute the integrals with respect to $\tau$. The following basic analytic fact will be useful, which we restate for convenience.

	\begin{lemma}\label{tau}
		Let $t>0$ and $0\leq p<q\leq 1$. Suppose the function $f:(pt,qt)\rightarrow\R$ satisfies the following two conditions:
		\begin{enumerate}
			\item There exists a unique $\tau_{max}$ such that $\tau_{max}=\argmax_{\tau\in[pt,qt]}f(\tau)$.
			\item For any $\tau\in[pt,qt]\setminus\{\tau_{max}\}$, $f'(\tau)(\tau_{max}-\tau)>0$.
		\end{enumerate}
		Then there exists some constant $C>0$ such that, as $t\rightarrow\infty$,
		\begin{equation}\label{tau_eq}
		\int_{pt}^{qt}\exp(f(\tau))\dd\tau=C\exp(f(\tau_{max})+o(t)).
		\end{equation}
	\end{lemma}
	

\section{proofs of the main theorems}\label{proof}
	
	\subsection{Proofs for the time-constrained probabilities.}\label{proof_time}
	
	In this subsection, we prove Theorem \ref{time} and \ref{time_late} using the estimates from Section 2. We first prove Theorem \ref{time}.
	
	\begin{proof}[Proof of Theorem \ref{time}]
		We rewrite \eqref{time_eq}, the equation to be proved, in the following three cases with different ranges of $\alpha$. 
		\begin{enumerate}[(i)]
		\item Given $\alpha\in[-\rho,1)$, as $t\rightarrow\infty$,
		\begin{equation}\label{time_alpha_large}
		\Prob{X_{max}(t)\leq \sqrt{2}\alpha t,\,\,X\in T_{[0,\gamma t]}}
		=\begin{cases}
		\eee^{-t\l{22}{\gamma}+o(t)},&\text{ if  }\,\gamma\in\left(0,\frac{1-\alpha}{2\sqrt{2}-1}\right],\\
		\eee^{-t\ll{11}{\gamma}{1}+o(t)},&\text{ if  }\,\gamma\in\left(\frac{1-\alpha}{2\sqrt{2}-1},\frac{1-\alpha}{\sqrt{2}}\right],\\
		\eee^{-t\ll{11}{\frac{1-\alpha}{\sqrt{2}}}{1}+o(t)},&\text{ if }\gamma\in\left(\frac{1-\alpha}{\sqrt{2}},1\right].
		\end{cases} 
		\end{equation}
		\item Given $\alpha\in(-2\rho,-\rho)$, as $t\rightarrow\infty$,
		\begin{equation}\label{time_alpha_middle}
		\Prob{X_{max}(t)\leq \sqrt{2}\alpha t,\,\,X\in T_{[0,\gamma t]}}
		=\begin{cases}
		\eee^{-t\l{32}{\gamma}+o(t)},&\text{ if  }\,\gamma\in\left(0,-\frac{\alpha+\rho}{\rho}\right],\\
		\eee^{-t\l{22}{\gamma}+o(t)},&\text{ if  }\,\gamma\in\left(-\frac{\alpha+\rho}{\rho},\frac{1-\alpha}{2\sqrt{2}-1}\right],\\
		\eee^{-t\ll{11}{\gamma}{1}+o(t)},&\text{ if }\gamma\in\left(\frac{1-\alpha}{2\sqrt{2}-1},1\right].
		\end{cases} 
		\end{equation}
		\item Given $\alpha\in(-\infty,-2\rho]$, as $t\rightarrow\infty$,
		\begin{equation}\label{time_alpha_small}
		\Prob{X_{max}(t)\leq \sqrt{2}\alpha t,\,\,X\in T_{[0,\gamma t]}}
		=\eee^{-t\l{32}{\gamma}+o(t)}\text{ for all }\gamma\in\left(0,1\right].
		\end{equation}
		\end{enumerate}
		
		Applying Lemma \ref{decomposition} with $A=[0,\gamma t]$ and $B=(-\infty,\infty)$, we rewrite the probability $\Prob{X_{max}(t)\leq \sqrt{2}\alpha t,\,\,X\in T_{[0,\gamma t]}}$ as 
		\begin{equation}\label{time_split_pre}
		\int_{0}^{\gamma t}Y_1\left(I_1\right)\dd\tau+\int_{0}^{\gamma t}Y_2\left(I_2\right)\dd\tau+\int_{0}^{\gamma t}Y_1\left(I_3\right)\dd\tau,
		\end{equation}
		which, after the change of variable $\lambda_\tau=\tau/t$, becomes
		\begin{equation}\label{time_split}
		\begin{split}
			&\int_{0}^{\gamma }tY_1\left(\left(-\infty,\sqrt{2}\alpha t-\sqrt{2}(1-\lambda_\tau)t\right]\right)\dd\lambda_\tau\\
			+\,\,&
			\int_{0}^{\gamma }tY_2\left(\sqrt{2}\alpha t+\left[-\sqrt{2}(1-\lambda_\tau)t,\sqrt{2}\rho(1-\lambda_\tau)t\right]\right)\dd\lambda_\tau\\
			+\,\,&\int_{0}^{\gamma }tY_3\left(\left[\sqrt{2}\alpha t+\sqrt{2}\rho(1-\lambda_\tau),\infty\right)\right)\dd\lambda_\tau.
		\end{split}
		\end{equation}
		In the remainder of the proof, we will estimate the three summands in \eqref{time_split}, and  then compare their orders.		
		Using Lemma \ref{y1} with $\lambda_u=1$, we can rewrite the first summand of \eqref{time_split} as
		\begin{equation}\label{time_y1_tau}
			\int_{0}^{\gamma\wedge\left(1-\alpha\right)}\mathrm{e}^{-t\ll{11}{\lambda_\tau}{1}+o(t)}\dd\lambda_\tau+\int_{\gamma\wedge\left(1-\alpha\right)}^{\gamma}\mathrm{e}^{-t\l{12}{\lambda_\tau} +o(t)}\dd\lambda_\tau.
		\end{equation}
		Observe that when regarded as a function of $\lambda_\tau$, $\ll{11}{\lambda_\tau}{1}$ strictly increases as $\lambda_\tau<\frac{1-\alpha}{\sqrt{2}}$ and strictly decreases as $\lambda_\tau>\frac{1-\alpha}{\sqrt{2}}$, while $\l{12}{\lambda_\tau}$ strictly decreases for all $\lambda_\tau$. By Lemma \ref{tau}, we know that \eqref{time_y1_tau} is equal to
		\begin{equation}
			\begin{cases}
			\mathrm{e}^{-t\ll{11}{\gamma}{1}+o(t)},&\text{ if }0<\gamma\leq\frac{1-\alpha}{\sqrt{2}},\\
			\mathrm{e}^{-t\ll{11}{\frac{1-\alpha}{\sqrt{2}}}{1}+o(t)}, &\text{ if } \frac{1-\alpha}{\sqrt{2}}<\gamma\leq 1-\alpha,\\
			\mathrm{e}^{-t\ll{11}{\frac{1-\alpha}{\sqrt{2}}}{1}+o(t)}+\mathrm{e}^{-t\l{12}{1-\alpha} +o(t)},&\text{ if }1-\alpha<\gamma\leq 1,
			\end{cases}
		\end{equation}
		which, as $-\ll{11}{\frac{1-\alpha}{\sqrt{2}}}{1}>-\l{12}{1-\alpha}$ since $2\rho<1$, is equal to
		\begin{equation}\label{time_1_result}
			\begin{cases}
			\mathrm{e}^{-t\ll{11}{\gamma}{1}+o(t)},&\text{ if }0<\gamma\leq\frac{1-\alpha}{\sqrt{2}},\\
			\mathrm{e}^{-t\ll{11}{\frac{1-\alpha}{\sqrt{2}}}{1}+o(t)}, &\text{ if } \frac{1-\alpha}{\sqrt{2}}<\gamma\leq 1.
			\end{cases}
		\end{equation}
		
		Using Lemma \ref{y2} with $\lambda_u=1$ and $\lambda_v=-\rho$, we rewrite the second summand of \eqref{time_split} as 
		\begin{equation}\label{time_y2_tau}
			\int_{0}^{\gamma\wedge\left(-\frac{\alpha+\rho}{\rho}\right)}\mathrm{e}^{-t\ll{21}{\lambda_\tau}{-\rho}+o(t)}\dd\lambda_\tau
			+\int_{\gamma\wedge\left(-\frac{\alpha+\rho}{\rho}\right)}^{\gamma\wedge\left(\frac{1-\alpha}{2\sqrt{2}-1}\right)}\mathrm{e}^{-t\l{22}{\lambda_\tau}+o(t)}\dd\lambda_\tau
			+\int_{\gamma\wedge\left(\frac{1-\alpha}{2\sqrt{2}-1}\right)}^{\gamma}\mathrm{e}^{-t\ll{21}{\lambda_\tau}{1}+o(t)}\dd\lambda_\tau.
		\end{equation}
		Notice that the first and the second exponents in \eqref{time_y2_tau}, regarded as functions of $\lambda_\tau$, both strictly increase for all $\lambda_\tau>0$. Since $\ll{21}{\cdot}{1}=\ll{11}{\cdot}{1}$, the monotonicity of the third exponent has been described below \eqref{time_y1_tau}. By Lemma \ref{tau}, \eqref{time_y2_tau} is equal to
		\begin{numcases}{}
			\mathrm{e}^{-t\ll{21}{\gamma}{-\rho}+o(t)}, & if $0<\gamma\leq-\frac{\alpha+\rho}{\rho}$,\label{time_y2_tau_1}\\
			\mathrm{e}^{-t\ll{21}{-\frac{\alpha+\rho}{\rho}}{-\rho}+o(t)}+\mathrm{e}^{-t\l{22}{\gamma}+o(t)}, & if $-\frac{\alpha+\rho}{\rho}<\gamma\leq\frac{1-\alpha}{2\sqrt{2}-1}$,\label{compare_BE'_C}\\
			\mathrm{e}^{-t\ll{21}{-\frac{\alpha+\rho}{\rho}}{-\rho}+o(t)}+\mathrm{e}^{-t\l{22}{\frac{1-\alpha}{2\sqrt{2}-1}}+o(t)}+\mathrm{e}^{-t\ll{21}{\gamma}{1}+o(t)}, & if $\frac{1-\alpha}{2\sqrt{2}-1}<\gamma\leq\frac{1-\alpha}{\sqrt{2}}$,\label{compare_A_BE'_C'}\\
			\mathrm{e}^{-t\ll{21}{-\frac{\alpha+\rho}{\rho}}{-\rho}+o(t)}+\mathrm{e}^{-t\l{22}{\frac{1-\alpha}{2\sqrt{2}-1}}+o(t)}+\mathrm{e}^{-t\ll{21}{\frac{1-\alpha}{\sqrt{2}}}{1}+o(t)},& if $\frac{1-\alpha}{\sqrt{2}}<\gamma\leq 1$.\label{compare_opt_BE'_C'}
		\end{numcases}
		Observe that the difference
		\begin{equation}
		-\ll{21}{-\frac{\alpha+\rho}{\rho}}{-\rho}-\left(-\l{22}{\gamma}\right)= -\left(4\sqrt{2}\rho-1\right)\gamma-4\sqrt{2}\left(\alpha+\rho\right)+\frac{\alpha+\rho}{\rho}
		\end{equation}
		is negative when $\gamma>-\frac{\alpha+\rho}{\rho}$, which implies that the second exponential terms in \eqref{compare_BE'_C}-\eqref{compare_opt_BE'_C'} are of larger order than the first terms. Moreover, the third term in \eqref{compare_A_BE'_C'} dominates the second, since
		\begin{equation}\label{compare_A_C'_diff}
		\begin{split}
			&-\ll{21}{\gamma}{1}-\left(-\l{22}{\frac{1-\alpha}{2\sqrt{2}-1}}\right)\\
			&\qquad=\,\,
			\frac{2}{\gamma}\left( \left(\gamma-\frac{11-4\sqrt{2}}{4(2\sqrt{2}-1)}\left(1-\alpha\right)\right)^2-\left(\left(\frac{11-4\sqrt{2}}{4(2\sqrt{2}-1)}\right)^2-\frac{1}{2} \right)\left(1-\alpha\right)^2 \right),
		\end{split}
		\end{equation}
		which is negative when
		\begin{equation}
		\frac{1-\alpha}{2\sqrt{2}-1}<\gamma<\left(\sqrt{2}-\frac{1}{2}\right)(1-\alpha),
		\end{equation}
		which is the case for $\gamma\in\left(\frac{1-\alpha}{2\sqrt{2}-1},\frac{1-\alpha}{\sqrt{2}} \right]$ in \eqref{compare_A_BE'_C'}. The third term in \eqref{compare_opt_BE'_C'} also dominates the second, since
		\begin{equation}
		-\ll{21}{\frac{1-\alpha}{\sqrt{2}}}{1}-\left(-\l{22}{\frac{1-\alpha}{2\sqrt{2}-1}}\right)=\frac{\rho^2}{2\sqrt{2}-1}(1-\alpha)>0.
		\end{equation}
		Hence, we conclude that the second term in \eqref{time_split} is of order
		\begin{equation}\label{time_2_result}
		\begin{cases}
		\mathrm{e}^{-t\ll{21}{\gamma}{-\rho}+o(t)}, &\text{ if }0<\gamma\leq-\frac{\alpha+\rho}{\rho},\\
		\mathrm{e}^{-t\l{22}{\gamma}+o(t)}, &\text{ if }-\frac{\alpha+\rho}{\rho}<\gamma\leq\frac{1-\alpha}{2\sqrt{2}-1},\\
		\mathrm{e}^{-t\ll{21}{\gamma}{1}+o(t)}, &\text{ if }\frac{1-\alpha}{2\sqrt{2}-1}<\gamma\leq\frac{1-\alpha}{\sqrt{2}},\\
		\mathrm{e}^{-t\ll{21}{\frac{1-\alpha}{\sqrt{2}}}{1}+o(t)}, &\text{ if }\frac{1-\alpha}{\sqrt{2}}<\gamma\leq 1.
		\end{cases}
		\end{equation}
		
		Applying Lemma \ref{y3} with $\lambda_v=-\rho$, we rewrite the third summand in \eqref{time_split} as
		\begin{equation}\label{time_y3_tau}
			\int_{0}^{\gamma\wedge\left(-\frac{\alpha+\rho}{\rho} \right)}\mathrm{e}^{-t\l{32}{\lambda_\tau}+o(t)}\dd\lambda_\tau+\int_{\gamma\wedge\left(-\frac{\alpha+\rho}{\rho}\right)}^{\gamma}\mathrm{e}^{-t\ll{31}{\lambda_\tau}{-\rho}+o(t)}\dd\lambda_\tau.
		\end{equation}
		Notice that the two exponents in \eqref{time_y3_tau}, regarded as functions of $\lambda_\tau$, strictly increase when $\lambda_\tau>0$. Thus by Lemma \ref{tau}, \eqref{time_y3_tau} is equal to
		\begin{numcases}{}
			\mathrm{e}^{-t\l{32}{\gamma}+o(t)}, & if $0<\gamma\leq-\frac{\alpha+\rho}{\rho}$,\\
			\mathrm{e}^{-t\l{32}{-\frac{\alpha+\rho}{\rho}}+o(t)}+\mathrm{e}^{-t\ll{31}{\gamma}{-\rho}+o(t)},& if $-\frac{\alpha+\rho}{\rho}<\gamma\leq 1$.\label{compare_BE_D'}
		\end{numcases}
		Notice that the second term in \eqref{compare_BE_D'} dominates, since
		\begin{equation}
		-\ll{31}{\gamma}{-\rho}-\left(-\l{32}{-\frac{\alpha+\rho}{\rho}}\right)=\frac{1+\rho^2}{\gamma}\left(\left(\gamma+\frac{\alpha+\rho}{2\rho(1+\rho^2)}\right)-\left(\frac{(2\rho^2+1)(\alpha+\rho)}{2\rho(1+\rho^2)}\right) \right),
		\end{equation}
		which is positive if
		\begin{equation}
		\gamma<\frac{\rho(\alpha+\rho)}{1+\rho^2}\,\text{  or  }\,\gamma>-\frac{\alpha+\rho}{\rho},
		\end{equation}
		satisfied by the range of $\gamma$ in \eqref{compare_BE_D'}. Thus the third term in \eqref{time_split} is of order
		\begin{equation}\label{time_3_result}
			\begin{cases}
			\mathrm{e}^{-t\l{32}{\gamma}+o(t)},&\text{ if }0<\gamma\leq-\frac{\alpha+\rho}{\rho},\\
			\mathrm{e}^{-t\ll{31}{\gamma}{-\rho}+o(t)}, &\text{ if }-\frac{\alpha+\rho}{\rho}<\gamma\leq 1.
			\end{cases}
		\end{equation}
		
		So far, we have obtained the estimates of the three summands in \eqref{time_split}, which are shown in \eqref{time_1_result}, \eqref{time_2_result}, and \eqref{time_3_result}. Next, we compare them for different ranges of $\alpha$ and $\gamma$. 
		
		\underline{Case (i).} Since $\alpha\in[-\rho,1)$, we have
		\begin{equation}
			-\frac{\alpha+\rho}{\rho}\leq 0<\frac{1-\alpha}{2\sqrt{2}-1}<\frac{1-\alpha}{\sqrt{2}}<1-\alpha<1.
		\end{equation}
		Adding \eqref{time_1_result}, \eqref{time_2_result}, and \eqref{time_3_result} together, \eqref{time_split} is equal to
		\begin{numcases}{}
			\mathrm{e}^{-t\ll{11}{\gamma}{1}+o(t)}+\mathrm{e}^{-t\l{22}{\gamma}+o(t)}+\mathrm{e}^{-t\ll{31}{\gamma}{-\rho}+o(t)}=\mathrm{e}^{-t\l{22}{\gamma}+o(t)},& if  $0<\gamma\leq\frac{1-\alpha}{2\sqrt{2}-1}$,\label{compare_A_BE_C}\\
			2\mathrm{e}^{-t\ll{11}{\gamma}{1}+o(t)}+\mathrm{e}^{-t\ll{31}{\gamma}{-\rho}+o(t)}=\mathrm{e}^{-t\ll{11}{\gamma}{1}+o(t)},& if  $\frac{1-\alpha}{2\sqrt{2}-1}<\gamma\leq\frac{1-\alpha}{\sqrt{2}}$,\label{compare_A_BE}\\
			2\mathrm{e}^{-t\ll{11}{\frac{1-\alpha}{\sqrt{2}}}{1}+o(t)}+\mathrm{e}^{-t\ll{31}{\gamma}{-\rho}+o(t)}=\mathrm{e}^{-t\ll{11}{\frac{1-\alpha}{\sqrt{2}}}{1}+o(t)},& if  $\frac{1-\alpha}{\sqrt{2}}<\gamma\leq 1$,\label{compare_opt_BE}
		\end{numcases}
		where the equality in \eqref{compare_A_BE_C} holds because 
		\begin{equation}
			-\l{22}{\gamma}-\left(-\ll{11}{\gamma}{1}\right)=\frac{9-4\sqrt{2}}{\gamma}\left(\gamma-\frac{1-\alpha}{2\sqrt{2}-1} \right)^2\geq 0\text{ for all }0<\gamma\leq 1,
		\end{equation}
		and
		\begin{equation}
			-\l{22}{\gamma}-\left(-\ll{31}{\gamma}{-\rho}\right)=\frac{1}{\gamma}(\rho\gamma+\alpha+\rho)^2>0\text{ for all }0<\gamma\leq 1,
		\end{equation}
		the equality in \eqref{compare_A_BE} holds because
		\begin{equation}\label{compare_A_BE_diff}
		-\ll{11}{\gamma}{1}-\left(-\ll{31}{\gamma}{-\rho}\right)=-\frac{6-2\sqrt{2}}{\gamma}\left(\left(\gamma-\frac{2-\sqrt{2}\alpha}{6-2\sqrt{2}} \right)^2-\frac{2(\alpha+2\rho)^2}{(6-2\sqrt{2})^2} \right),
		\end{equation}
		which is nonnegative when
		\begin{equation}
		\frac{2-2\sqrt{2}(\alpha+\rho)}{6-2\sqrt{2}}\leq\gamma\leq 1,
		\end{equation}
		satisfied by the corresponding $\gamma$ range as $\frac{2-2\sqrt{2}(\alpha+\rho)}{6-2\sqrt{2}}<\frac{1-\alpha}{2\sqrt{2}-1}$ for $\alpha>-2\rho$, and the equality in \eqref{compare_opt_BE} holds by \eqref{compare_A_BE} and
		\begin{equation}
			\mathrm{e}^{-t\ll{11}{\frac{1-\alpha}{\sqrt{2}}}{1}+o(t)}\geq \mathrm{e}^{-t\ll{11}{\gamma}{1}+o(t)}.
		\end{equation}
	 	Then, \eqref{time_alpha_large} follows directly from \eqref{compare_A_BE_C}-\eqref{compare_opt_BE}.
		
		\underline{Case (ii).} Since $\alpha\in(-2\rho,-\rho)$, we have
		\begin{equation}
			0<-\frac{\alpha+\rho}{\rho}<\frac{1-\alpha}{2\sqrt{2}-1}<1<\frac{1-\alpha}{\sqrt{2}}<1-\alpha.
		\end{equation}
		Thus \eqref{time_split} is equal to the sum of \eqref{time_1_result}, \eqref{time_2_result}, and \eqref{time_3_result}, 
		\begin{numcases}{}
			\mathrm{e}^{-t\ll{11}{\gamma}{1}+o(t)}+\mathrm{e}^{-t\ll{21}{\gamma}{-\rho}+o(t)}+\mathrm{e}^{-t\l{32}{\gamma}+o(t)}=\mathrm{e}^{-t\l{32}{\gamma}+o(t)},& if $0<\gamma<-\frac{\alpha+\rho}{\rho}$,\label{compare_A_BE_D}\\
			\mathrm{e}^{-t\ll{11}{\gamma}{1}+o(t)}+\mathrm{e}^{-t\l{22}{\gamma}+o(t)}+\mathrm{e}^{-t\ll{31}{\gamma}{-\rho}+o(t)}=\mathrm{e}^{-t\l{22}{\gamma}+o(t)},& if $-\frac{\alpha+\rho}{\rho}\leq\gamma\leq\frac{1-\alpha}{2\sqrt{2}-1}$,\label{time_alpha_middle_gamma_middle}\\
			2\mathrm{e}^{-t\ll{11}{\gamma}{1}+o(t)}+\mathrm{e}^{-t\ll{21}{\gamma}{-\rho}+o(t)}=\mathrm{e}^{-t\ll{11}{\gamma}{1}+o(t)},& if $\frac{1-\alpha}{2\sqrt{2}-1}<\gamma\leq 1$,\label{time_alpha_middle_gamma_large}
		\end{numcases}
		where the equality in \eqref{time_alpha_middle_gamma_middle} follows from \eqref{compare_A_BE_C}, the equality in \eqref{time_alpha_middle_gamma_large} follows from \eqref{compare_A_BE}, and the equality in \eqref{compare_A_BE_D} is because of the following two facts. The third term in \eqref{compare_A_BE_D} is of larger order than the second term, since
		\begin{equation}\label{compare_BE_D_diff}
		-\l{32}{\gamma}-\left(-\ll{21}{\gamma}{-\rho}\right)=\frac{1-\gamma}{\gamma(1+\gamma)}\left(\rho\gamma+(\alpha+\rho)\right)^2\geq 0\text{ for all }0<\gamma\leq 1.
		\end{equation}
		In addition, the second term in \eqref{compare_A_BE_D} dominates the first term, since $\ll{21}{\cdot}{-\rho}=\ll{31}{\cdot}{-\rho}$ and we learn from \eqref{compare_A_BE_diff} that
		\begin{equation}
		-\ll{21}{\gamma}{-\rho}-\left(-\ll{11}{\gamma}{1}\right)\leq 0 \iff \frac{2-2\sqrt{2}(\alpha+\rho)}{6-2\sqrt{2}}\leq\gamma\leq 1,
		\end{equation}
		and $-\frac{\alpha+\rho}{\rho}<\frac{2-2\sqrt{2}(\alpha+\rho)}{6-2\sqrt{2}}$ when $\alpha>-2\rho$. Combining \eqref{compare_A_BE_D}-\eqref{time_alpha_middle_gamma_large}, \eqref{time_alpha_middle} follows.
		
		\underline{Case (iii).} Since $\alpha\in(-\infty,-2\rho]$, we have 
		\begin{equation}
			0<1\leq\frac{1-\alpha}{2\sqrt{2}-1}\leq-\frac{\alpha+\rho}{\rho}.
		\end{equation}
		Thus for all $0<\gamma\leq 1$, \eqref{time_split} is equal to the sum of \eqref{time_1_result}, \eqref{time_2_result}, and \eqref{time_3_result},
		\begin{equation}\label{compare_A_BE_D_alpha_small}
			\mathrm{e}^{-t\ll{11}{\gamma}{1}+o(t)}+\mathrm{e}^{-t\ll{21}{\gamma}{-\rho}+o(t)}+\mathrm{e}^{-t\l{32}{\gamma}+o(t)}=\mathrm{e}^{-t\l{32}{\gamma}+o(t)},
		\end{equation}
		where the equality holds by the following two facts. Firstly, the first and the second exponents in \eqref{compare_A_BE_D_alpha_small} have appeared in \eqref{compare_A_BE} with different $\alpha$ ranges. The difference between these two exponents is recorded in \eqref{compare_A_BE_diff}, which, under the condition that $\alpha<-2\rho$, is not positive as that would require
		 \begin{equation}
		 1<\gamma<\frac{2-2\sqrt{2}(\alpha+\rho)}{6-2\sqrt{2}}.
		 \end{equation}
		Thus the second term in \eqref{compare_A_BE_D_alpha_small} is of larger order than the first term. Secondly, from \eqref{compare_BE_D_diff}, we know that the third term in \eqref{compare_A_BE_D_alpha_small} dominates the second one for all $0<\gamma\leq 1$. Hence, \eqref{compare_A_BE_D_alpha_small} holds and thus also \eqref{time_alpha_small} follows.  This concludes the proof of Theorem \ref{time}.
	\end{proof}
	Next, we prove Theorem \ref{time_late}, which focuses on the $\alpha\in(-\rho,1)$ and restricts the first branching time $\tau$ to be later than the optimal one $\frac{1-\alpha}{\sqrt{2}}t$.
	
	\begin{proof}[Proof of Theorem \ref{time_late}]
		In the same way of obtaining \eqref{time_split_pre}-\eqref{time_split} in the proof of Theorem \ref{time}, we apply Lemma \ref{decomposition} with $A=[\gamma t, t]$ and $B=(-\infty,\infty)$ and the change of variables $\lambda_\tau=\tau/t$ to rewrite $\Prob{X_{max}(t)\leq \sqrt{2}\alpha t,\,\,X\in T_{[\gamma t,t]}}$ as
		\begin{equation}\label{time_late_split}
		\begin{split}
		&\int_{\gamma }^{1}tY_1\left(\left(-\infty,\sqrt{2}\alpha t-\sqrt{2}(1-\lambda_\tau)t\right]\right)\dd\lambda_\tau\\
		+\,\,&\int_{\gamma }^{1}tY_2\left(\left[\sqrt{2}\alpha t-\sqrt{2}(1-\lambda_\tau)t,\sqrt{2}\alpha t+\sqrt{2}\rho(1-\lambda_\tau)t\right]\right)\dd\lambda_\tau\\
		+\,\,&\int_{\gamma }^{1}tY_3\left(\left[\sqrt{2}\alpha t+\sqrt{2}\rho(1-\lambda_\tau)t,\infty\right)\right)\dd\lambda_\tau.
		\end{split}
		\end{equation}
		
		Notice that we restrict $-\rho<\alpha<1$ and $\frac{1-\alpha}{\sqrt{2}}<\gamma\leq 1$ in the theorem. For the first summand in \eqref{time_late_split}, we apply Lemma \ref{y1} with $\lambda_u=1$ to rewrite it as
		\begin{equation}\label{time_late_y1_lemma}
			\int_{\gamma}^{\gamma\vee \left(1-\alpha\right)}\mathrm{e}^{-t\ll{11}{\lambda_\tau}{1}+o(t)}\dd\lambda_\tau+\int_{\gamma\vee \left(1-\alpha\right)}^{1}\mathrm{e}^{-t\l{12}{\lambda_\tau}+o(t)}\dd\lambda_\tau.
		\end{equation}
		By Lemma \ref{tau}, \eqref{time_late_y1_lemma} is equal to
		\begin{equation}
			\begin{cases}
				\mathrm{e}^{-t\ll{11}{\gamma}{1}+o(t)}+\mathrm{e}^{-t\l{12}{1-\alpha}+o(t)},&\text{ if }0<\gamma<1-\alpha,\\
				\mathrm{e}^{-t\l{12}{\gamma}+o(t)}, &\text{ if }1-\alpha\leq\gamma\leq 1.
			\end{cases}
		\end{equation}
		Observe that
		\begin{equation}
			-\ll{11}{\gamma}{1}-\left(-\l{12}{1-\alpha} \right)= -\frac{2}{\gamma}\left(\left(\gamma-\frac{3(1-\alpha)}{4}\right)^2-\left(\frac{1-\alpha}{4}\right)^2 \right),
		\end{equation}
		which is positive when $\frac{1-\alpha}{2}<\gamma<1-\alpha$, satisfied by our $\gamma$ range. Thus \eqref{time_late_y1_lemma} is equal to
		\begin{equation}\label{time_late_1_result}
		\begin{cases}
			\mathrm{e}^{-t\ll{11}{\gamma}{1}+o(t)},&\text{ if }0<\gamma<1-\alpha,\\
			\mathrm{e}^{-t\l{12}{\gamma}+o(t)}, &\text{ if }1-\alpha\leq\gamma\leq 1.
		\end{cases}
		\end{equation}
		
		For the second summand in \eqref{time_late_split}, since $\lambda_\tau\geq\gamma>\frac{1-\alpha}{2\sqrt{2}-1}$, by Lemma \ref{y2} with $\lambda_u=1$, $\lambda_v=-\rho$ and Lemma \ref{tau}, this summand is equal to
		\begin{equation}\label{time_late_2_result}
			\int_{\gamma }^{1}\mathrm{e}^{-t\ll{21}{\lambda_\tau}{1}+o(t)}\dd\lambda_\tau=\mathrm{e}^{-t\ll{21}{\gamma}{1}+o(t)}.
		\end{equation}
		
		For the third summand in \eqref{time_late_split}, since $\lambda_\tau\geq\gamma>-\frac{\alpha+\rho}{\rho}$, by Lemma \ref{y3} with $\lambda_v=-\rho$ and Lemma \ref{tau}, this summand is equal to
		\begin{equation}\label{time_late_3_result}
			\int_{\gamma }^{1}\mathrm{e}^{-t\ll{31}{\lambda_\tau}{-\rho}+o(t)}\dd\lambda_\tau=\mathrm{e}^{-t\ll{31}{1}{-\rho}+o(t)}.
		\end{equation}
		
		To conclude the proof, we show that \eqref{time_late_1_result} dominates \eqref{time_late_2_result} and \eqref{time_late_3_result}. Since $\ll{11}{\cdot}{1}=\ll{21}{\cdot}{1}$, it suffices to show that $-\l{12}{\gamma}\geq-\ll{21}{\gamma}{1}$ and $-\ll{21}{\gamma}{1}\geq-\ll{31}{1}{-\rho}$. The former is obviously true for all $0<\gamma\leq 1$. For the latter, 
		\begin{equation}
		-\ll{21}{\gamma}{1}-\left(-\ll{31}{1}{-\rho}\right)=-\frac{2}{\gamma}\left(\left(\gamma-\frac{\alpha^2-2\alpha+3}{4}\right)^2-\left(\frac{\alpha^2-2\alpha+3}{4}\right)^2+\frac{(1-\alpha)^2}{2} \right),
		\end{equation}
		which is nonnegative when
		\begin{equation}
		\frac{(1-\alpha)^2}{2}\leq \gamma\leq 1.
		\end{equation}
		As our range of $\gamma$ is a subset of this, our claim is verified. Hence \eqref{time_late_1_result} indeed dominates over \eqref{time_late_2_result} and \eqref{time_late_3_result} and the proof is done.
	\end{proof}

	\subsection{Proofs for the location-constrained probabilities.}\label{proof_location}
	
	In this subsection, we prove Theorem \ref{location_below} and Theorem \ref{location_above}.
	We first prove Theorem \eqref{location_below}, which restricts the first branching location to be below $\sqrt{2}\alpha t-\sqrt{2}\beta(t-\tau)$, for some $\beta\geq 1$.
	
	\begin{proof}[Proof of Theorem \ref{location_below}]
		Let $A=\left[(\gamma-\epsilon)t,\gamma t\right]$ and $B=\left(-\infty,\sqrt{2}\alpha t-\sqrt{2}\beta(t-\tau)\right]$. Since $\beta\geq 1$, 
			$\sqrt{2}\alpha t-\sqrt{2}\beta (t-\tau)\leq\sqrt{2}\alpha t-\sqrt{2}(t-\tau)$
		and thus with the notation from Lemma \ref{decomposition}
		\begin{equation}\label{loc.b.1}
			B\cap I_1=B, \,\,B\cap I_2=B\cap I_3=\emptyset.
		\end{equation}
		Using Lemma \ref{decomposition},  we rewrite $\Prob{X_{max}(t)\leq \sqrt{2}\alpha t,\,X\in T_{[(\gamma-\epsilon)t,\gamma t]}\cap L_{\left(-\infty,\sqrt{2}\alpha t-\sqrt{2}\beta(t-\tau)\right]}}$ as
		\begin{equation}\label{location_below_tau}
			\int_{(\gamma-\epsilon)t}^{\gamma t}Y_1\left(\left(-\infty,\sqrt{2}\alpha t-\sqrt{2}\beta(t-\tau)\right]\right)\dd\tau = 	\int_{\gamma-\epsilon}^{\gamma }tY_1\left(\left(-\infty,\sqrt{2}\alpha t-\sqrt{2}\beta(1-\lambda_\tau)t\right]\right)\dd\lambda_\tau
,		 		\end{equation}
		after the change of variables $\lambda_\tau=\tau/t$.
		Applying Lemma \ref{y1} with $\lambda_u=\beta$, we obtain that
		\begin{equation}
			Y_1\left(\left(-\infty,\sqrt{2}\alpha t-\sqrt{2}\beta(1-\lambda_\tau)t\right]\right)=\begin{cases}
			\mathrm{e}^{-t\ll{11}{\lambda_\tau}{\beta}+o(t)},&\text{ if }\beta>\beta_1^{\alpha}(\lambda_\tau),\\
			\mathrm{e}^{-t\l{12}{\lambda_\tau}+o(t)},&\text{ if }\beta\leq\beta_1^{\alpha}(\lambda_\tau).
			\end{cases}
		\end{equation}
		Thus, by Lemma \ref{tau}, \eqref{location_below_tau} is equal to
		\begin{equation}\label{time_late_result_pre}
			\begin{cases}
			\mathrm{e}^{-t\ll{11}{\gamma-\epsilon}{\beta}+o(t)},&\text{ if }\beta>\beta_1^{\alpha}(\gamma),\\
			\mathrm{e}^{-t\l{12}{\gamma-\epsilon}+o(t)},&\text{ if }\beta\leq\beta_1^{\alpha}(\gamma).
			\end{cases}
		\end{equation}
		The desired results follow if we take $t\rightarrow\infty$ and then $\epsilon\rightarrow 0$ in \eqref{time_late_result_pre}.
	\end{proof}
	
	Next, we prove Theorem \ref{location_above}, where the first branching position is constrained to be above $\sqrt{2}\alpha t-\sqrt{2}\beta(t-\tau)$, for some $\beta\leq 1$.
	
	\begin{proof}[Proof of Theorem \ref{location_above}]
		Let $A=[(\gamma-\epsilon)t,\gamma t]$ and $B=\left[\sqrt{2}\alpha t-\sqrt{2}\beta(t-\tau),\infty\right)$. From $\beta\leq 1$, we know that $B\cap I_1=\emptyset$. To obtain the values of $B\cap I_2$ and $B\cap I_3$, we need to divide the $\beta$ range into two and discuss two cases.
		
		\underline{Case 1: $\beta\in\left(-\infty,-\rho\right]$.} In this $\beta$ range, we have that
		\begin{equation}
			B\cap I_1=\emptyset, \,\,B\cap I_2=\emptyset, \,\,B\cap I_3=B.
		\end{equation}		
		Then by Lemma \ref{decomposition}, $\Prob{X_{max}(t)\leq \sqrt{2}\alpha t,\,X\in T_{[(\gamma-\epsilon)t,\gamma t]}\cap L_{\left[\sqrt{2}\alpha t-\sqrt{2}\beta(t-\tau),\infty\right)}}$ can be rewritten as
		\begin{equation}\label{location_above_beta_small_tau}
			\int_{(\gamma-\epsilon)t}^{\gamma t}Y_3\left(\left[\sqrt{2}\alpha t-\sqrt{2}\beta(t-\tau),\infty\right) \right)\dd\tau =\int_{\gamma-\epsilon}^{\gamma }tY_3\left(\left[\sqrt{2}\alpha t-\sqrt{2}\beta(1-\lambda_\tau)t,\infty\right) \right)\dd\lambda_\tau,
		\end{equation}
		after changing variables $\lambda_\tau=\tau/t$. 
		Applying Lemma \ref{y3} with $\lambda_v=\beta$ and Lemma \ref{tau}, \eqref{location_above_beta_small_tau} is equal to
		\begin{equation}
			\begin{cases}
			\mathrm{e}^{-t\ll{31}{\gamma}{\beta}+o(t) },&\text{ if }\beta<\frac{\alpha}{1+\gamma},\\
			\mathrm{e}^{-t\l{32}{\gamma}+o(t)},&\text{ if }\beta\geq\frac{\alpha}{1+\gamma}.
			\end{cases}
		\end{equation}
		As
			$\frac{\alpha}{1+\gamma}<-\rho\iff\gamma<-\frac{\alpha+\rho}{\rho}$,
		we summarize the results in Case 1 as follows.
		\begin{itemize}
			\item If $0<\gamma<-\frac{\alpha+\rho}{\rho}$, then the estimate for the desired probability is of order
			\begin{equation}\label{location_above_beta_small_gamma_small}
				\begin{cases}
				\mathrm{e}^{-t\ll{31}{\gamma}{\beta}+o(t) },&\text{ if }\beta<\frac{\alpha}{1+\gamma},\\
				\mathrm{e}^{-t\l{32}{\gamma}+o(t)},&\text{ if }\frac{\alpha}{1+\gamma}\leq\beta\leq 1.
				\end{cases}
			\end{equation}
			
			\item If $-\frac{\alpha+\rho}{\rho}\leq\gamma<1$, then the estimate is of order
			\begin{equation}\label{location_above_beta_small_gamma_large}
				\mathrm{e}^{-t\ll{31}{\gamma}{\beta}+o(t) },\text{ for all }\beta\leq 1.
			\end{equation}
		\end{itemize}
	
		\underline{Case 2: $\beta\in\left(-\rho,1\right]$.} In this $\beta$ range,
		\begin{equation}
			B\cap I_2=\left[\sqrt{2}\alpha t-\sqrt{2}\beta(t-\tau),\sqrt{2}\alpha t+\sqrt{2}\rho(t-\tau)\right], \,\,B\cap I_3=I_3.
		\end{equation}
		Thus by Lemma \ref{decomposition}, $\Prob{X_{max}(t)\leq \sqrt{2}\alpha t,\,X\in T_{[(\gamma-\epsilon)t,\gamma t]}\cap L_{\left[\sqrt{2}\alpha t-\sqrt{2}\beta(t-\tau),\infty\right)}}$ can be rewritten as
		\begin{equation}
		\begin{split}
			&\int_{(\gamma-\epsilon)t}^{\gamma t}Y_2\left(\left[\sqrt{2}\alpha t-\sqrt{2}\beta(t-\tau),\sqrt{2}\alpha t+\sqrt{2}\rho(t-\tau)\right) \right)\dd\tau\\
			+\,\,&\int_{(\gamma-\epsilon)t}^{\gamma t}Y_3\left(\left[\sqrt{2}\alpha t+\sqrt{2}\rho(t-\tau),\infty\right) \right)\dd\tau,
		\end{split}
		\end{equation}
		which, after the change of variables, becomes
		\begin{equation}\label{location_above_beta_large_tau}
		\begin{split}
			&\int_{\gamma-\epsilon}^{\gamma }tY_2\left(\left[\sqrt{2}\alpha t-\sqrt{2}\beta(1-\lambda_\tau)t,\sqrt{2}\alpha t+\sqrt{2}\rho(1+\lambda_\tau)t\right) \right)\dd\lambda_\tau\\
			+\,\,&\int_{\gamma-\epsilon}^{\gamma }tY_3\left(\left[\sqrt{2}\alpha t+\sqrt{2}\rho(1-\lambda_\tau)t,\infty\right) \right)\dd\lambda_\tau.
		\end{split}
		\end{equation}
		
		Applying Lemma \ref{y2} with $\lambda_u=\beta$ and $\lambda_v=-\rho$ and Lemma \ref{tau}, the first summand in \eqref{location_above_beta_large_tau} is equal to
		\begin{equation}\label{location_above_beta_large_y2}
			\begin{cases}
			\mathrm{e}^{-t\ll{21}{\gamma}{-\rho}+o(t) },&\text{ if }0<\gamma<-\frac{\alpha+\rho}{\rho},\\
			\mathrm{e}^{-t\l{22}{\gamma}+o(t) },&\text{ if }-\frac{\alpha+\rho}{\rho}\leq\gamma\leq\frac{\beta-\alpha}{\beta+2\rho},\\
			\mathrm{e}^{-t\ll{21}{\gamma}{\beta}+o(t)},&\text{ if }\frac{\beta-\alpha}{\beta+2\rho}<\gamma<1.
			\end{cases}
		\end{equation}
		
		By Lemma \ref{y3} with $\lambda_v=-\rho$ and Lemma \ref{tau}, the second summand in \eqref{location_above_beta_large_tau} is equal to
		\begin{equation}\label{location_above_beta_large_y3}
			\begin{cases}
				\mathrm{e}^{-t\l{32}{\gamma}+o(t)}, &\text{ if }0<\gamma<-\frac{\alpha+\rho}{\rho},\\
				\mathrm{e}^{-t\ll{31}{\gamma}{-\rho}+o(t) },&\text{ if }-\frac{\alpha+\rho}{\rho}\leq\gamma<1.
			\end{cases}
		\end{equation}
		Note that we put the equality sign in the different conditions compared to Lemma \ref{y3}, since the two cases are equal when $\lambda_\tau=-\frac{\alpha+\rho}{\rho}$. 
		
		Adding together \eqref{location_above_beta_large_y2} and \eqref{location_above_beta_large_y3}, we obtain that \eqref{location_above_beta_large_tau} is equal to
		\begin{numcases}{}
			\mathrm{e}^{-t\ll{21}{\gamma}{-\rho}+o(t) }+\mathrm{e}^{-t\l{32}{\gamma}+o(t)}=\mathrm{e}^{-t\l{32}{\gamma}+o(t)},& if $0<\gamma<-\frac{\alpha+\rho}{\rho}$,\label{compare'_BE_D}\\
			\mathrm{e}^{-t\l{22}{\gamma}+o(t) }+\mathrm{e}^{-t\ll{31}{\gamma}{-\rho}+o(t) }=\mathrm{e}^{-t\l{22}{\gamma}+o(t) },& if $-\frac{\alpha+\rho}{\rho}\leq\gamma\leq\frac{\beta-\alpha}{\beta+2\rho}$,\label{compare'_BE_C}\\
			\mathrm{e}^{-t\ll{21}{\gamma}{\beta}+o(t)}+\mathrm{e}^{-t\ll{31}{\gamma}{-\rho}+o(t) }=\mathrm{e}^{-t\ll{21}{\gamma}{\beta}+o(t)},& if $\frac{\beta-\alpha}{\beta+2\rho}<\gamma<1$,\label{compare_beta_BE}
		\end{numcases}
		where the equalities in \eqref{compare'_BE_D} and \eqref{compare'_BE_C} follow from \eqref{compare_A_BE_D} and \eqref{compare_A_BE_C}, respectively. The equality in \eqref{compare_beta_BE} holds since
		\begin{equation}
		\begin{split}
		&-\ll{21}{\gamma}{\beta}-\left(-\ll{31}{\gamma}{-\rho}\right)\\
		=\,&-\frac{(1-\gamma)^2}{\gamma}\left(\beta-\frac{2\rho\gamma+\alpha}{1-\gamma} \right)^2+\left(3-2\sqrt{2} \right)\gamma+\frac{(\alpha+\rho)^2}{\gamma}+2\rho\alpha-4\rho+2\\
		\geq\,&-\ll{21}{\gamma}{-\rho}-\left(-\ll{31}{\gamma}{-\rho}\right)=0,
		\end{split}
		\end{equation}
		where the inequality is due to the facts that $-\rho<\beta\leq 1$ and, if $\gamma\neq 1$,
		\begin{equation}
			\gamma>\frac{\beta-\alpha}{\beta+2\rho}\iff \beta<\frac{\alpha+2\rho\gamma}{1-\gamma}.
		\end{equation} 
		Hence, \eqref{compare'_BE_D}-\eqref{compare_beta_BE} follow and we obtain the desired estimate for Case 2.
		Grouping the estimates in \eqref{location_above_beta_small_gamma_small}-\eqref{location_above_beta_small_gamma_large} for Case 1  and \eqref{compare'_BE_D}-\eqref{compare_beta_BE} for Case 2 using the fact that
		\begin{equation}
			\gamma<-\frac{\alpha+\rho}{\rho}\iff\frac{\alpha}{1+\gamma}<-\rho
		\end{equation}
		and, when $\beta>-\rho$,
		\begin{equation}
			\alpha>-2\rho\iff -\frac{\alpha+\rho}{\rho}<\frac{\beta-\alpha}{\beta+2\rho}<1,
		\end{equation}
		we obtain all estimates as stated in the theorem.
	\end{proof}

	\textbf{Acknowledgements. } We would like to express our gratitude to Anton Bovier for his continued support and valuable suggestions throughout this project.

\end{document}